\documentclass[a4paper,10pt]{article}
\usepackage[latin1]{inputenc}
\usepackage[T1]{fontenc}                  
\usepackage[english]{babel}
\usepackage{verbatim}
\usepackage{amsmath}
\usepackage{amsthm, amsfonts, amssymb}
\usepackage{amssymb}
\usepackage{graphicx}
\usepackage{bm}
\usepackage{fancyhdr}                  
\usepackage{color}

\theoremstyle{definition}
\newtheorem{defin}{Definition}[section]

\newtheorem{rem}[defin]{Remark}

\theoremstyle{plane}
\newtheorem{thm}[defin]{Theorem}
\newtheorem{prop}[defin]{Proposition}
\newtheorem{coroll}[defin]{Corollary}
\newtheorem{lemma}[defin]{Lemma}

\newcommand{\mbb}{\mathbb}

\newcommand{\mc}{\mathcal}
\newcommand{\veps}{\varepsilon}
\newcommand{\vrho}{\varrho}

\newcommand{\wtilde}{\widetilde}
\newcommand{\vphi}{\varphi}
\newcommand{\oline}{\overline}
\newcommand{\ra}{\rightarrow}

\newcommand{\hra}{\hookrightarrow}

\newcommand{\eps}{\varepsilon}
\newcommand{\R}{\mathbb{R}}

\newcommand{\N}{\mathbb{N}}
\newcommand{\Z}{\mathbb{Z}}

\renewcommand{\div}{{\rm div}\,}

\newcommand{\Id}{{\rm Id}\,}

\newcommand{\Int}{\displaystyle \int}

\def\d{\partial}

\def\dj{\Delta_j}
\def\tilde{\widetilde}
\def\hat{\widehat}
\def\div{{\rm div}\,}

\def\cF{{\mathcal F}}

\def\cP{{\mathcal P}}
\def\cQ{{\mathcal Q}}
\def\cR{{\mathcal R}}
\def\cS{{\mathcal S}}

\renewcommand{\div}{{\rm div}\,}

\def\d{\partial}
\def\div{{\rm div}\,}

\usepackage{enumerate}
\usepackage{xcolor}

\def\d{\partial}

\def\dj{\Delta_j}

\def\tilde{\widetilde}
\def\hat{\widehat}
\def\div{{\rm div}\,}

\def\cF{{\mathcal F}}

\def\cP{{\mathcal P}}
\def\cQ{{\mathcal Q}}
\def\cS{{\mathcal S}}

\allowdisplaybreaks

\author{
Francesco FANELLI 
 and
 Xian LIAO\thanks{Corresponding author.}}

\date\today

\begin{document}

\title{The well-posedness issue  for an   inviscid zero-Mach number system in general Besov spaces}

\maketitle

\subsubsection*{Abstract}
{\small The present paper is devoted to the study of a zero-Mach number system with heat conduction but no viscosity.
We work in the framework of general non-homogeneous Besov spaces $B^s_{p,r}(\R^d)$, with $p\in[2,4]$ and for any  $d\geq 2$, which can be
embedded into the class of globally Lipschitz functions.

We prove a local in time well-posedness result in these classes for general initial densities and velocity fields.
Moreover, we are able to show  a continuation criterion and a lower bound for the lifespan of the solutions.

The proof of the results relies on Littlewood-Paley decomposition and paradifferential calculus, and on refined commutator estimates
in Chemin-Lerner spaces.}

\medbreak

\noindent\textbf{Keywords.}
{\small Zero-Mach number system; well-posedness; Besov spaces; Chemin-Lerner spaces; continuation criterion; lifespan.}

\medbreak
\noindent\textbf{Mathematics Subject Classification (2010).}
{\small Primary: 35Q35. Secondary: 76N10, 35B65.}

\section{Introduction}

The free evolution  of a compressible,   effectively  heat-conducting but inviscid fluid  obeys the following equations:
\begin{equation}\label{eq:fullnewtonian}\left\{
\begin{array}{ccc}
\d_t\rho+\div(\rho v)&=&0,\\
\d_t(\rho v)+\div(\rho v\otimes v)+\nabla p&=&0,\\
\d_t(\rho e)+\div(\rho v e)-\div(k\nabla\vartheta)+p\,\div v&=&0,
\end{array}
\right.
\end{equation}
where $\rho=\rho(t,x)\in\R^+$
 stands for the mass density, $v=v(t,x)\in\R^d$ for  the velocity field and $e=e(t,x)\in\R^+$
for  the internal energy per unit mass. The time variable $t$
belongs to $\R^+$ or to $[0,T]$ and the space variable $x$ is in $\R^d$
with $d\geq2.$
The scalar functions  $p=p(t,x)$ and $\vartheta=\vartheta(t,x)$ denote the pressure and temperature respectively.
The heat-conducting coefficient $k=k(\rho,\vartheta)$ is supposed to be smooth in both its variables.

We supplement   System \eqref{eq:fullnewtonian} with the following two state equations:
\begin{equation*}\label{eq:state}
p=R\rho \vartheta,\quad e=C_v\vartheta,
\end{equation*}
where $R,C_v$ denote the ideal gas constant and the specific heat capacity at constant volume, respectively. That is, we restrict ourselves to (so-called) ideal gases.


In this paper, we will consider highly subsonic ideal gases strictly away from vacuum, and correspondingly, we will work
with the \emph{inviscid zero-Mach number system} (see   \eqref{eq:EH} or  \eqref{system}  below) which derives from System \eqref{eq:fullnewtonian} by letting
the Mach number go to zero.
For completeness, in the following we derive the zero-Mach number system \emph{formally}.

Just as in \cite{D-L}, suppose  $(\rho,v,p)$  to be a solution of System \eqref{eq:fullnewtonian} and  define the dimensionless Mach number $\eps$ to be the ratio
of the velocity $v$ by the reference sound speed.  Then  the rescaled triplet
$$
\left(\rho_\eps(t,x)=\rho\Bigl(\frac{t}{\eps},x\Bigr),\quad
v_\eps(t,x)=\frac{1}{\eps}v\Bigl(\frac{t}{\eps},x\Bigr),\quad
p_\eps(t,x)=p\Bigl(\frac{t}{\eps},x\Bigr)\right)
$$
satisfies the following non-dimensional system
\begin{equation}\label{eq:epsilon}\left\{
\begin{array}{ccc}
\d_t\rho_\eps+\div(\rho_\eps v_\eps)&=&0,\\[1ex]
\d_t(\rho_\eps v_\eps)+\div (\rho_\eps v_\eps\otimes v_\eps)+\frac{\nabla p_\eps}{\eps^2}&=&0,\\[1ex]
\frac{1}{\gamma-1}(\d_tp_\eps+\div(p_\eps v_\eps ))-\div(k_\eps\nabla \vartheta_\eps)+p_\eps\div v_\eps&=&0.
\end{array}
\right.
\end{equation}
Here   $\gamma:=C_p/{C_v}=1+R/{C_v}$ represents the adiabatic index and the constant $C_p$ denotes the specific heat capacity at constant pressure.
 The rescaled temperature and heat-conducting coefficient are given by
$$
\vartheta_\eps(t,x)=\vartheta\Bigl(\frac{t}{\eps},x\Bigr),\quad k_\eps(t,x)=\frac{1}{\eps}k\Bigl(\frac{t}{\eps},x\Bigr).
$$
Now let $\eps$ go to $0$, that is, the pressure $p_\eps$ equals to a positive constant $P_0$ by Equations $\eqref{eq:epsilon}_2$
 and $\eqref{eq:epsilon}_3$: thus System \eqref{eq:epsilon} becomes formally the following zero Mach number system
immediately (see \cite{Alazard06}, \cite{PLions96} and \cite{Zeytounian04} for detailed computations):
\begin{equation}\label{eq:EH}
\left\{ \begin{array}{ccc}
\d_t \rho+\div(\rho v) & =&0,\\
 \d_{t}(\rho v)+\div(\rho v\otimes v ) + \nabla\Pi & = & 0, \\
 \div v - \frac{\gamma -1}{\gamma P_0}\div( k \nabla\vartheta)& = &0,
\end{array} \right.
\end{equation}
 where  $\Pi=\Pi(t,x)$ is a new unknown function.

Let us  make reference to  some results on  the \textit{incompressible limit} of Euler equations: namely, incompressible Euler equations
can be viewed as compressible Euler equations when the Mach number tends to vanish.
For the \textit{isentropic} Euler system, there are many early works such as  \cite{Ebin77,  Isozaki87, K-M, Schochet, Ukai}, etc.
At the beginning of this century,  \cite{Alazard05, M-S} treated   the \textit{non-isentropic} case.
Later  Alazard \cite{Alazard06}  generalized the study to  various models, which include the case of the low Mach number limit from
System \eqref{eq:fullnewtonian} to System \eqref{eq:EH}.

 \smallbreak

Next, we will reformulate System \eqref{eq:EH} into a new system (see System \eqref{system} below) with a new \emph{divergence-free} velocity field.
We will mainly deal with this new system in this paper.
 Before going into details, let us suppose here that the
density $\rho$ always has positive lower bound and
 converges to some constant (say ``1'') at infinity, in the sense detailed below (see \eqref{initial data}). We also mention
that in the following, the fact that $\rho\vartheta\equiv P_0/R$ is a positive constant   will be used thoroughly.

Similarly as in \cite{D-L},   for notational simplicity, we set $\alpha$ to be the positive constant defined by
$$
\alpha\,=\,\frac{\gamma-1}{\gamma P_0}\,=\,\frac{R}{C_p P_0}\,=\,\frac{1}{C_p \rho\vartheta}.
$$
Then we define the following two coefficients, always viewed as regular functions of $\rho$ only:
\begin{equation}\label{kappa,lambda}
\kappa=\kappa(\rho)=\alpha k\vartheta\quad \hbox{and}\quad \lambda=\lambda(\rho)=\rho^{-1}.
\end{equation}
One furthermore introduces two scalar functions $a=a(\rho)$ and $b=b(\rho)$,  such that
\begin{equation}\label{relation:a,b}
\nabla a\,=\,\kappa\nabla\rho
\,=\,-\rho\nabla b,\quad a(1)=b(1)=0.
\end{equation}

We then define the new ``velocity'' $u$ and the new ``pressure'' $\pi$  respectively as
\begin{equation}\label{relation:u}
u\,=\,v-\alpha k\nabla\vartheta
  \,\equiv\,v-\nabla b
  \,\equiv\,v+\kappa\rho^{-1}\nabla\rho,
\qquad \pi \,=\,\Pi-\d_t a\,.
\end{equation}
Then System \eqref{eq:EH} finally becomes
\begin{equation}\label{system}
\left\{
\begin{array}{cc}
&\d_t\rho+u\cdot\nabla\rho-\div(\kappa\nabla\rho)=0,\\
&\d_t u+(u+\nabla b)\cdot\nabla u+\lambda\nabla\pi=h,\\
&\div u=0,
\end{array}
\right.
\end{equation}
where $\kappa, b, \lambda$ are defined above (see \eqref{kappa,lambda} and \eqref{relation:a,b}) and
\begin{equation}\label{h}
h(\rho,u)\,=\,\rho^{-1}\div (v\otimes \nabla a).
\end{equation}
Let us just verify Equation $\eqref{system}_2$.
 Observe that $(\ref{eq:EH})_2$ gives
\begin{equation}\label{eq:origin,u}
\d_t(\rho u)+\div(\rho v\otimes u)+\d_t( \rho\nabla b)+\div(v\otimes  \rho\nabla b)+\nabla\Pi=0.
\end{equation}
It is easy to find that
$$\d_t(\rho\nabla b)+\nabla\Pi=-\d_t\nabla a+\nabla\Pi=\nabla \pi\quad\hbox{and}\quad
  \div(v\otimes \rho\nabla b)=-\div(v\otimes \nabla a).$$
Thus by view of  Equation $\eqref{eq:EH}_1$, Equation \eqref{eq:origin,u} can be rewritten as
\begin{equation}\label{eq:u}
\rho\d_t u+\rho v\cdot\nabla u+\nabla\pi=\div(v\otimes\nabla a).
\end{equation}
 We thus multiply \eqref{eq:u} by $\lambda=\rho^{-1}$   to get Equation $\eqref{system}_2$.

\smallbreak

In System \eqref{system}, although the ``velocity'' $u$ is divergence-free, one encounters  a (quasilinear) parabolic equation for the density $\rho$ and  the ``source'' term $h$   involves two derivatives of $\rho$.
Note that if simply $\kappa\equiv 0$, then $a\equiv b\equiv 0$ and hence System \eqref{system} becomes the so-called
density-dependent Euler equations
\begin{equation}\label{system:dens-depend}
\left\{
\begin{array}{cc}
&\d_t\rho+v\cdot\nabla\rho =0,\\
&\d_t v+v \cdot\nabla v+\lambda\nabla\pi=0,\\
&\div v=0.
\end{array}
\right.
\end{equation}
In the above \eqref{system:dens-depend}, due to the null heat conduction,  $\rho$ satisfies a transport equation,  $h$  vanishes and the velocity $v$ itself is solenoidal.
As early as in 1980,  Beir\~{a}o da Veiga and  Valli   \cite{B-V1,B-V2}  investigated  \eqref{system:dens-depend}.
We also cite the book \cite{Antontsev-K-M} as a good survey of the boundary-value problems for nonhomogeneous fluids.
By use of an energy identity,   Danchin  \cite{D} studied System \eqref{system:dens-depend} in the framework of nonhomogeneous Besov space
$B^s_{p,r}(\R^d)$ which can be embedded in $C^{0,1}$.
Recently in \cite{D-F}, Danchin and the first author treated the end point case $B^s_{\infty,r}$, and studied the lifespan of the solutions
in the case of space dimension $d=2$.

If the fluid is viscous, that is to say there is an additional viscous stress tensor  within the momentum equation $\eqref{eq:EH}_2$,
then System \eqref{eq:EH} becomes the low Mach number limit system of the full Navier-Stokes system.
See   \cite{Alazard06, D-L, L}  and references therein for some relevant results.

However, to our knowledge, there are few well-posedness results for the  inviscid low Mach-number limit system  \eqref{eq:EH}.
Notice that System  \eqref{eq:EH} can also be viewed as a \emph{nonhomogeneous system in the presence of diffusion}, which describes
an inviscid  fluid consisting of two components (say, water and salt), both incompressible, with a mass diffusion effect between them
(the so-called Fick's law):
$$
\div (v+\kappa\nabla\ln\rho)=0.
$$
In this case,
 $\rho$ and $v$ are considered to be the mean density and the \textit{mean-mass} velocity of the mixture respectively,
 $\kappa$ denotes the positive diffusion coefficient,
 and $\nabla\Pi$, as usual, denotes some unknown pressure.
For more physical backgrounds of this model, see \cite{F-K}.
One can also refer to  Beir\~ao da Veiga \textit{et al.} \cite{B-S-V} for an existence-uniqueness result of \textit{classical} solutions.

\medbreak
In this paper, we will study the well-posedness of the Cauchy problem for System \eqref{eq:EH} in the
framework of \textit{general Besov spaces} $B^s_{p,r}(\R^d)$ (with $p\in[2,4]$ and in any space dimension $d\geq 2$)
which can be embedded into the class of globally Lipschitz functions.

Similarly as in \cite{D,D-F}, the analysis will be based on an intensive use of the para-differential calculus
and some  (newly-developed) commutator estimates.

Moreover, refined a priori estimates lead us to state a continuation criterion (in the same spirit of the well-known
result of \cite{Beale-K-M84} for the homogeneous incompressible Euler equations), and to find a lower bound for the lifespan of the solutions
in terms of the initial data only.

We refer to Section \ref{s:results} for more details on our working hypothesis and  discussions on the obtained results.
Let us just say here that the restriction $p\in[2,4]$, which is due to the analysis of the pressure term $\nabla\pi$, can be somehow
relaxed.
For instance, for \emph{finite-energy} initial data $(\rho_0,u_0)$,
well-posedness for system \eqref{system} can be recovered for any \emph{$1<p\leq+\infty$}.
We refer to \cite{F-L_L2} for an analysis in this direction: the endpoint case $B^1_{\infty,1}$ is permitted there and in dimension $d=2$  the lower bound for the lifespan
is refined such that the solutions tend to be globally defined for initial densities which are small perturbations of a constant state.


\medbreak
Our paper is organized in the following way.

In  next section we will present our main  local-in-time well-posedness result Theorem  \ref{th:w-p}.
We will also state a continuation criterion   and a lower bound for the lifespan  in Theorems \ref{th:cc_r} and \ref{th:lifespan} respectively.

Section \ref{s:tools} is devoted to the tools  from Fourier analysis.

In Section \ref{s:th1} we will tackle the proof of Theorem \ref{th:w-p}: we will give some fundamental commutator
estimates (see Lemma \ref{l:tilde}) and product estimates (Lemma \ref{l:prod}) in the \textit{time-dependent} Besov spaces.

Section \ref{s:conti-life}  is devoted to the proof of Theorem \ref{th:cc_r}
and  Theorem \ref{th:lifespan}.

Finally, in the appendix we will give the complete  proof of   Lemma \ref{l:tilde}.

\section{Main results} \label{s:results}
Let us focus on System \eqref{system} to introduce our main results. In view of Equation $\eqref{system}_1$, of parabolic type,
by  \textit{maximum principle} we can assume that the density $\rho$ (if it exists on the time interval $[0,T]$)  has the same positive
 upper and lower bounds as the initial density $\rho_0$:
\begin{equation*}\label{bound}
0<\rho_\ast\leq \rho(t,x) \leq \rho^\ast,\quad \forall\, t\in [0,T],\, x\in \R^d.
\end{equation*}
Correspondingly, the coefficients $\kappa$ and $\lambda$ can always be bounded from above and below, which ensures that
the pressure $\pi$ satisfies an \textit{elliptic} equation in divergence form: applying operator ``$\div$'' to Equation $\eqref{system}_2$
gives the following:
\begin{equation}\label{eq:pi}
\div(\lambda\,\nabla\pi)\,=\,\div(h\,-\,v\cdot\nabla u)\,,\qquad \lambda\geq \lambda_\ast\,=\,(\rho^\ast)^{-1}\,>\,0\,.
\end{equation}


For system  \eqref{system}  there is no gain of regularity for the velocity $u$: we then suppose the initial divergence-free
``velocity'' field $u_0$ to belong to some space
$B^s_{p,r}$ which can be continuously embedded in $C^{0,1}$, i.e. the triplet $(s,p,r)\in \R\times [1,+\infty]^2$ has to satisfy
the following condition:
\begin{equation}\label{index:s,r}
s\,>\,1\,+\,\frac{d}{p}\;,\qquad\qquad\mbox{or}\qquad\qquad s\,=\,1\,+\,\frac{d}{p}\;,\;\;r\,=\,1\,.
\end{equation}
This  requires the ``source'' term $h$ to belong to $ L^1([0,T];B^s_{p,r})$
which, by view of definition \eqref{h} of $h$,  asks at least
$\nabla^2 \rho \in L^1([0,T];B^s_{p,r})$.
Keeping in mind that $\rho$ satisfies the parabolic type  equation, we expect to   gain   two orders of regularity (in space) when
taking the average in time.
We thus   have to assume the initial inhomogeneity $\rho_0-1$  to be in the same space $B^s_{p,r}$ as above.
 However in general, we only get $\nabla^2\rho$ in the time-dependent Besov space $\tilde L^1_T(B^s_{p,r})$, which is a little bit larger
  than $L^1_T(B^s_{p,r})$ (see Definition \ref{def:Besov,tilde}).
Therefore in the whole paper we will deal rather with the spaces $\tilde L^\infty_T(B^s_{p,r})$
and $\tilde L^1_T(B^s_{p,r})$ (first introduce in \cite{Chemin-Lerner} by Chemin and Lerner);
in particular,  in Section  \ref{s:th1} we will give new commutator estimates and product estimates in these
\textit{time-dependent} Besov spaces, which imply a priori estimates for System \eqref{system}.

On the other hand, in order to control the \textit{low frequencies}  for $\nabla\pi$, one has to make sure that $h-v\cdot \nabla u\in L^2$.
Indeed, for Equation \eqref{eq:pi} above, the a priori estimate
$$\lambda_\ast\|\nabla\pi\|_{L^q}\leq C\|h-v\cdot\nabla u\|_{L^q}$$
holds \textit{independently} of $\lambda$
  only when $q=2$ (see Lemma 2 of \cite{D}).
     Hence the fact that
$h$ is composed of quadratic forms entails that $p$ has to verify
\begin{equation}\label{index:p}
p\in [2,4].
\end{equation}

To conclude, we have the following theorem, whose proof will be shown in Section \ref{s:th1}. 
\begin{thm} \label{th:w-p}
 Let the triplet $(s,p,r)\in\R^3$ satisfy conditions \eqref{index:s,r} and \eqref{index:p}.
Let us take an initial density state $\rho_0$ and an initial   velocity field $u_0$ such that
\begin{equation}\label{initial data}
0<\rho_\ast\leq \rho_0\leq \rho^\ast\,,\quad \div u_0=0,\quad \|\rho_0-1\|_{B^s_{p,r}}\,+\,\|u_0\|_{B^s_{p,r}}\,\leq\,M\,,
\end{equation}
 for some positive constants $\rho_\ast$, $\rho^\ast$ and $M$.
Then there exist a positive time $T$ (depending only on $\rho_\ast, \rho^\ast, M,d,s,p,r$) and a unique solution $(\rho,u,\nabla\pi)$ to System \eqref{system} such
that $(\vrho,u,\nabla\pi):=(\rho-1,u,\nabla\pi)$  belongs
to the space $E^s_{p,r}(T)$,  defined as the set of  triplet $ (\varrho, u, \nabla \pi) $ such that
\begin{equation}\label{space:E}
\left\{\begin{array}{c }
\varrho \in \wtilde C([0,T];B^{s}_{p,r})\cap\wtilde L^{1}([0,T];B^{s+2}_{p,r})\,,
\quad \rho_\ast\leq \varrho+1\leq\rho^\ast, \\[1ex]
u   \in   \wtilde{C}([0,T];B^{s}_{p,r})\,, \\[1ex]
\nabla\pi   \in  \wtilde{L}^{1}([0,T];B^{s}_{p,r})\cap L^1([0,T];L^2)\,,
\end{array}\right.
\end{equation}
with $\tilde C_{w}([0,T];B^s_{p,r}) \hbox{ if }\,r=+\infty $ (see also Definition \ref{def:Besov,tilde}).
\end{thm}

\begin{rem}\label{rem:origin}
Let us state briefly here the corresponding well-posedness result for the original system \eqref{eq:EH}. By view of the change of
variables \eqref{relation:u}, we have $u=\cP v$, $\nabla b=\cQ v$, where $\cP$ denotes the Leray projector over divergence-free vector
fields and $\cQ=\Id -\cP$:
$\hat{\cQ u}(\xi) = -  (\xi/{|\xi|^2})\,\xi\cdot\hat u(\xi).$
Assume Conditions \eqref{index:s,r} and \eqref{index:p}, and  the initial datum $(\rho_0, v_0)$ such that
$$
0<\rho_\ast\leq \rho_0\leq \rho^\ast\,,\quad \nabla b(\rho_0)=\cQ v_0,\qquad \|\rho_0-1\|_{B^s_{p,r}}\,+\,\|\cP v_0\|_{B^s_{p,r}}\,\leq\,M\,.
$$
Then,  there exist a positive time $T$   and a unique solution $(\rho,v,\nabla\Pi)$ to System \eqref{eq:EH} such that
\begin{equation*}
\left\{\begin{array}{c}
\rho-1 \in  \wtilde C([0,T];B^{s}_{p,r})\cap\wtilde L^{1}([0,T];B^{s+2}_{p,r})\,, \\[1ex]
v  \in  \wtilde{C}([0,T];B^{s-1}_{p,r})\cap \tilde L^2([0,T];B^s_{p,r})\,, \quad \cP v\in \tilde {C}([0,T];B^s_{p,r})\, ,\\[1ex]
\nabla\Pi  \in  \wtilde{L}^{1}([0,T];B^{s}_{p,r})\,,
\end{array}\right.
\end{equation*}
with $\tilde C_{w}([0,T];B^s_{p,r}) \hbox{ if }\,r=+\infty $.

One notices from above that the initial velocity $v_0$ needs not to be in $B^s_{p,r}$ whereas the velocity $v(t)$  will belong to it for almost every $t\in [0,T]$.
 On the other side,   this local existence result would hold  if, initially, $\rho_0-1\in B^{s+1}_{p,r}$, $v_0\in B^s_{p,r}$,   $\rho_0\in [\rho_\ast,\rho^\ast]$ and $\nabla b(\rho_0)=\cQ v_0$.

   Let us also point out that since $\d_t\rho\not\in L^1([0,T];L^2)$ in general, we do not know whether $\nabla\Pi\in L^1([0,T];L^2)$
(recall also definition \eqref{relation:u}).
    Hence it seems not convenient to deal with System \eqref{eq:EH} directly since the low frequences of $\nabla\Pi$ can not be controlled \textit{a priori}.
\end{rem}
\begin{rem} \label{r:p}
If we assume an additional \textit{smallness} hypothesis over the initial inhomogeneity aside from \eqref{initial data}, which ensures that
the pressure satisfies a Laplace equation (up to a perturbation term), then Condition \eqref{index:p} imposed on $p$ is  not necessary.
Theorem \ref{th:w-p} still holds true,  except for the fact  $\nabla\pi\in L^1([0,T];L^2)$.
\end{rem}


Next,  one can get a Beale-Kato-Majda type continuation criterion (see \cite{Beale-K-M84} for the original version)
for solutions to System \eqref{system}.
Notice that, according to  the solutions space $E^s_{p,r}$ defined by \eqref{space:E}, in order to  bound the $\tilde L^1_T(B^s_{p,r})$-norm
of the \textit{nonlinear} terms in $  h $ and
$v\cdot\nabla u $, one requires   $ (u,\nabla\rho)\in {L^\infty_T(L^\infty)}$, $(\nabla u, \nabla^2 \rho)\in {L^2_T(L^\infty)}$, etc.
Hence, by a refined a priori estimate (see Subsection \ref{s:conti}) we get the following statement.

\begin{thm}{\bf{[Continuation Criterion]}} \label{th:cc_r}
 Let the triplet $(s,p,r)\in\R\times[1,+\infty]^2$ satisfy conditions \eqref{index:s,r} and \eqref{index:p}. Let $(\rho,u,\nabla\pi)$ be a solution of \eqref{system} on $[0,T[\,\times\R^d$
such that:
$$
(\rho-1, u, \nabla\pi)\in \Bigl( \wtilde{C}([0,T[\,;B^s_{p,r})\cap\wtilde{L}^1_{loc}([0,T[\,;B^{s+2}_{p,r})\Bigr)
\times
\wtilde{C}([0,T[\,;B^s_{p,r})
\times
\wtilde{L}^1_{loc}([0,T[\,;B^s_{p,r})
$$
and, for some $\sigma>0$,
$$
 \sup_{t\in[0,T[}\|(\nabla\rho(t), u(t))\|_{L^\infty}\,+\,\int^T_0
\left(\left\|(\nabla^2\rho, \nabla u)\right\|^2_{L^\infty} +\|\nabla\pi(t)\|_{B^{-\sigma}_{p,\infty}\cap L^\infty}\right)\,dt\,<\,+\infty\,.
$$
Then $(\rho,u,\nabla\Pi)$ could be continued beyond $T$ (if $T$ is finite) into a solution of
\eqref{system} with the same regularity.
\end{thm}


Even in the two-dimensional case,
it's hard to expect global-in-time well-posedness for this system: the parabolic equation $\eqref{system}_1$
allows to improve regularity for the density term, but such a gain is (roughly speaking) deleted by the nonlinear
term in the momentum equation $\eqref{system}_2$.
 However, similar as in  \cite{D-F}, we manage to establish an \textit{explicit} lower bound for the lifespan of the solution, in \textit{any} dimension $d\geq2$. The proof will be the matter of Section \ref{s:life}.
\begin{thm} \label{th:lifespan}
 Under the hypotheses of Theorem \ref{th:w-p}, there exist  positive constants $L$, $\ell>6$ (depending only on $d$, $p$,
$r$, $\rho_*$ and $\rho^*$)
such that  the lifespan  of the solution to System \eqref{system} given by Theorem \ref{th:w-p} is bounded from below by the quantity
\begin{equation} \label{est:lifespan}
\frac{L}{1\,+\,\|u_0\|_{B^s_{p,r}}\,+\,\|\varrho_0\|^\ell_{B^s_{p,r}}}\,.
\end{equation}
\end{thm}

\begin{rem}\label{rem:cc_r}
Thanks to Theorem \ref{th:cc_r}, the lifespan is independent of the regularity.
Therefore, the $ {B^s_{p,r}}$-norm in \eqref{est:lifespan} can be replaced by the (weaker) ${B^{1+d/4}_{4,1}}$ norm.
\end{rem}

\begin{rem}
The lower bound \eqref{est:lifespan} might be improved by scaling.
Indeed, System  \eqref{system} is invariant under the following transformation:
\begin{equation*} \label{def:sol}
 (\rho^{ \veps }\,,\,u^{ \veps }\,,\,\nabla\pi^{ \veps } )(t\,,\,x)\,:=\,
 (\rho\,,\,\veps^{-1}\,u\,,\,\veps^{-2}\,\nabla\pi )(\veps^{-2} \,t\,, \eps^{-1}\,x).
\end{equation*}
Set $\eps^2= \|u_0\|_{B^{1+d/4}_{4,1}}$, then, by \eqref{est:lifespan}, the lifespan of the solution
$( \rho^\eps,  u^\eps,\nabla \pi^\eps)$ is bigger (up to a constant factor) than the quantity
$(1+\eps^{-1}\| \rho_0-1\|_{B^{1+d/4}_{4,1}}^{\ell})^{-1}$.
 Correspondingly the lifespan of $(\rho, u, \nabla\pi)$ is bigger than or equal to
$$\frac{\eps^{-2}L}{
 1+\eps^{-1}\| \rho_0-1\|_{B^{1+d/4}_{4,1}}^{\ell}
 }\,.$$
\end{rem}

\medbreak
We change the word here that in the sequel, $C$ always denotes some ``harmless'' constant (may vary from time to time)
depending only on $d,s,p,r,\rho_\ast,\rho^\ast$,
unless otherwise defined. Notation $A\lesssim B$ means $A\leq C B$ and $A\sim B$ says $A$ equals to $B$, up to a constant factor.
For notational convenience, the notation $\vrho$ always represent $\rho-1$, unless otherwise specified.

\section{An overview on Fourier analysis techniques} \label{s:tools}

Our results mostly rely on  Fourier analysis methods
which are based on  a nonhomogeneous dyadic partition of unity
with respect to   Fourier variable, the so-called Littlewood-Paley decomposition.
Unless otherwise specified, all the results which are presented
in this section are proved in \cite{B-C-D}, Chap. 2.

In order to define a  Littlewood-Paley decomposition,
fix a smooth radial function
$\chi$ supported in (say) the ball $B(0,\frac43),$
equals to $1$ in a neighborhood of $B(0,\frac34)$
and such that  $\chi$ is nonincreasing
over $\R_+$.
Set
$\varphi(\xi)=\chi(\frac\xi2)-\chi(\xi).$
{}
The {\it dyadic blocks} $(\Delta_j)_{j\in\Z}$
 are defined by\footnote{Throughout we agree  that  $f(D)$ stands for
the pseudo-differential operator $u\mapsto\cF^{-1}(f(\xi) \cF u(\xi)).$}
$$
\dj:=0\ \hbox{ if }\ j\leq-2,\quad\Delta_{-1}:=\chi(D)\quad\hbox{and}\quad
\Delta_j:=\varphi(2^{-j}D)\ \text{ if }\  j\geq0.
$$
We  also introduce the following low frequency cut-off:
$$
S_ju:=\chi(2^{-j}D)=\sum_{j'\leq j-1}\Delta_{j'}\quad\text{for}\quad j\geq0,\quad S_j u\equiv 0 \quad\text{for}\quad j\leq0.
$$
One can now define what a Besov space $B^s_{p,r}$ is:
\begin{defin}
\label{def:besov}
  Let  $u$ be a tempered distribution, $s$ a real number, and
$1\leq p,r\leq\infty.$ We set
$$
\|u\|_{B^s_{p,r}}:=\bigg(\sum_{j} 2^{rjs}
\|\Delta_j  u\|^r_{L^p}\bigg)^{\frac{1}{r}}\ \text{ if }\ r<\infty
\quad\text{and}\quad
\|u\|_{B^s_{p,\infty}}:=\sup_{j}\left( 2^{js}
\|\Delta_j  u\|_{L^p}\right).
$$
We then define the space $B^s_{p,r}$ as  the
subset of  distributions $u\in {\cS}'$ such  that
$\|u\|_{B^s_{p,r}}$ is finite.
\end{defin}

When solving evolutionary PDEs, it is natural to use
spaces of type $L^\rho_T(X)=L^\rho(0,T;X)$ with $X$ denoting some Banach space. In our case, $X$ will be a  Besov space so that we will
have to localize the equations by Littlewood-Paley decomposition. This will provide us
with  estimates of the Lebesgue norm of  each dyadic block \emph{before}
performing integration in time. This  leads to the following  definition for the so-called Chemin-Lerner Spaces, introduced for the first time in paper \cite{Chemin-Lerner}.
\begin{defin}\label{def:Besov,tilde}
For  $s\in \R$, $(q,p,r)\in [1,+\infty]^3$ and $T\in [0,+\infty]$, we set
$$
\|u\|_{\tilde L^q_T(B^s_{p,r})}
=\Bigl\| \Bigl(  2^{js}  \|\dj u(t)\|_{L^q_T(L^p)} \Bigr)_{j\geq -1}\Bigr\|_{\ell^r}.
$$
We also set $\tilde C_T(B^s_{p,r})=\tilde L_T^\infty(B^s_{p,r})\cap C([0,T];B^s_{p,r}).$
\end{defin}

\begin{rem} \label{rem:besov}
From the above definition, it is easy to see   $B^s_{2,2}\equiv H^s$ for all $s\in \R$ and $B^s_{p,r}\hookrightarrow C^{0,1}$ under
the hypothesis \eqref{index:s,r}. More generally, one has the continuous embedding
$B^k_{p,1}\hookrightarrow W^{k,p}\hookrightarrow B^k_{p,\infty}$ and a time-dependent version
 \begin{align*}\label{embed}
 L^1_t(B^{s}_{p,1})\,\hookleftarrow\,\tilde L^1_t(B^{s+\epsilon/2}_{p,1})\,\hookleftarrow\,
\tilde L^1_t(B^{s+\epsilon}_{p,\infty}),\,\forall \epsilon>0.
\end{align*}
 \end{rem}

  The following fundamental lemma (referred in what follows as \emph{Bernstein's inequalities})
  describes the way derivatives act on spectrally localized functions.
  \begin{lemma}\label{lpfond}
Let  $0<r<R.$
There exists a constant $C$  such that, for any $k\in \Z^+$, $\lambda\in\R^+$,  $(p,q)\in [1,\infty]^2$ with  $p\leq q$
and any function $u$ of~$L^p$, one has
$$
\displaylines{
{\rm Supp}\, \widehat u \subset   B(0,\lambda R)
\Longrightarrow
\|\nabla^k u\|_{L^q} \leq
 C^{k+1}\lambda^{k+d(\frac{1}{p}-\frac{1}{q})}\|u\|_{L^p};\cr
{\rm Supp}\, \widehat u \subset \{\xi\in\R^N\,/\, r\lambda\leq|\xi|\leq R\lambda\}
\Longrightarrow C^{-k-1}\lambda^k\|u\|_{L^p}
\leq
\|\nabla^k u\|_{L^p}
\leq
C^{k+1}  \lambda^k\|u\|_{L^p}.
}$$
\end{lemma}
  Lemma \ref{lpfond}  implies the following embedding result immediately, as a generalization of Remark \ref{rem:besov}:
  \begin{coroll}\label{c:embed}
  Space $B^{s_1}_{p_1,r_1}$ is continuously embedded in Space $B^{s_2}_{p_2,r_2}$ whenever
  $1\leq p_1\leq p_2\leq\infty$ and
  $$
  s_2< s_1-d/p_1+d/p_2\quad\hbox{or}\quad
  s_2=s_1-d/p_1+d/p_2\ \hbox{ and }\ 1\leq r_1\leq r_2\leq\infty.
  $$
  \end{coroll}

 Let us now recall the so-called Bony's decomposition introduced in \cite{Bony} for the products.
Formally, any product  of two tempered distributions $u$ and $v,$ may be decomposed
into
\begin{equation}\label{eq:bony}
uv=T_uv+T_vu+R(u,v)
\end{equation}
with
$$
T_uv:=\sum_j S_{j-1}u\dj v\ \hbox{ and }\
R(u,v):=\sum_j\sum_{|j'-j|\leq1}\dj u\,\Delta_{j'}v.
$$
The above operator $T$ is called ``paraproduct'' operator whereas
$R$ is called ``remainder'' operator.

We have the following classical estimates in Besov spaces
for the  products:
\begin{prop}\label{p:prod}
Let  $s,s_1,s_2\in \R$, $1\leq r,r_1,r_2,p \leq\infty$  with
 $\frac{1}{r}\leq\min\{1,\frac{1}{r_1}+\frac{1}{r_2}\}$.
\begin{itemize}
\item
For the paraproduct  one has
\begin{eqnarray*} 
\|T_u v\|_{B^{s}_{p,r}}\lesssim
\|u\|_{L^{\infty}}\|v\|_{B^{s}_{p,r}}\cdotp
\end{eqnarray*}
\item For the remainder, if $s_1+s_2+d\min\{0,1-\frac{2}{p}\}>0$, then
$$
\|R(u,v)\|_{B^{s_1+s_2- \frac{d}{p}}_{p,r}}\lesssim \|u\|_{B^{s_1}_{p,r_1}}\|v\|_{B^{s_2}_{p,r_2}}.
$$
\item If $s>0$, one has the following product estimate
$$
\|u v\|_{\tilde L^q_T(B^s_{p,r})}
\leq C\left(\|u\|_{L^{q_1}_T(L^\infty)}\|  v\|_{\tilde L^{q_2}_T(B^{s}_{p,r})}
\,+\, \|u\|_{\tilde L^{q_3}_T(B^{s}_{p,r})} \|  v\|_{L^{q_4}_T(L^\infty)}\right),
\quad  \frac 1q:=\frac{1}{q_1}  +\frac{1}{q_2}=\frac{1}{q_3} +\frac{1}{q_4}.
$$
\end{itemize}
\end{prop}

One also needs the following commutator estimate for the transport terms in the equations:
\begin{prop}\label{p:comm}
Let $s>-d\min\{1/p,1/{p'}\}$, $(p,r)\in [1,\infty]^2$ such that $r=1$ if $s=1+d/p$. Then we have
 \begin{equation*}
 \int^t_0\left\|2^{j s}\left\|[\vphi,\Delta_j]\nabla\psi\right\|_{L^p}\right\|_{\ell^r}d\tau\;\leq\;
 C\Int^t_0 \|\nabla \varphi(t)\|_{B^{\frac dp}_{p,1}\cap B^{s-1}_{p,r}}\|\nabla\psi\|_{B^{s-1}_{p,r}}\, d\tau.
 \end{equation*}
\end{prop}

The following results pertain to the composition of functions in Besov spaces: they will be needed for estimating functions
depending on the density. We refer to Chap. 2 of \cite{B-C-D} and to \cite{D} for their proofs.
\begin{prop}\label{p:comp_grad}
Let $I$ be an open  interval of $\R$ and $F:I\rightarrow\R$ a smooth function.
Then for any compact subset $J\subset I,$ $s>0$, $(q,p,r)\in[1,+\infty]^3$ and any function $a$ valued in $J$, we have
$$
\|\nabla(F(a))\|_{\tilde L^q_T(B^{s-1}_{p,r})}\leq C\|\nabla a\|_{\tilde L^q_T(B^{s-1}_{p,r})}.
$$
If furthermore   $F(0)=0$, then
$   \left\|F(a)\right\|_{\tilde L^q_T(B^s_{p,r})}\,\leq\,C\,\|a\|_{\tilde L^q_T(B^s_{p,r})}\,.
$
\end{prop}

Finally, we shall make an extensive use of energy  estimates (see \cite{D} for the proof)  for the
following elliptic equation satisfied by the pressure $\Pi$:
\begin{equation*}\label{eq:elliptic}
-\div(\lambda\nabla\Pi)=\div F\quad\hbox{in }\ \R^d,
\quad \lambda=\lambda(x)\geq \lambda_\ast>0,
\end{equation*}
\begin{lemma}\label{l:laxmilgram}
For all vector field $F$ with coefficients in $L^2,$ there exists a tempered distribution $\Pi,$
unique up to  constant functions,
such that
\begin{equation*}\label{eq:el0}
\lambda_\ast \|\nabla\Pi\|_{L^2}\leq\|F\|_{L^2}.
\end{equation*}
\end{lemma}



\section{Proof of the local well-posedness}\label{s:th1}


In  this section we aim to prove the well-posedness result  Theorem \ref{th:w-p}  for System \eqref{system}.
The a priori estimates we will establish in Subsection \ref{subs:est} will be the basis throughout the following context.

\subsection{Linearized equations}\label{subs:est}

In this subsection we will build a priori estimates for the linearized equations associated to  System \eqref{system}.
Firstly, one gives some  commutator and product estimates in the \textit{time-dependent} Besov spaces, which entails the a priori estimate for the density immediately.
\begin{lemma} \label{l:tilde}
Let $s, s_i, \sigma_i\in \R$, $i=1,2$, $(p,r)\in[1,+\infty]^2$.
\begin{itemize}
\item If $s>0$, then
\begin{equation}\label{linest:comm}
\int^t_0\left\|2^{j s}\left\|[\vphi,\Delta_j]\nabla\psi\right\|_{L^p}\right\|_{\ell^r}d\tau\;\leq\;
C\Int^t_0\left(\left\|\nabla\vphi\right\|_{L^\infty}\,\left\|\psi\right\|_{B^{s}_{p,r}}\,+\,
\left\|\nabla\vphi\right\|_{B^{s-1}_{p,r}}\,\left\|\nabla\psi\right\|_{L^\infty}\right)d\tau.
\end{equation}
\item If $s>0$ and
 \begin{equation*}
s+1=\theta s_1+(1-\theta)s_2=\eta\sigma_1+(1-\eta)\sigma_2,
\,\hbox{ with }\theta, \eta\in (0,1],
 \end{equation*}
then for any $\eps>0$, one has
\begin{align}
 &  \Bigl\|2^{js}\int^t_0
\left\|\nabla\left(\left[\vphi,\Delta_j\right]\nabla\psi\right)\right\|_{L^p}\,d\tau
\Bigr\|_{\ell^r}
\,\leq \,\frac{C\,\theta}{\veps^{(1-\theta)/\theta}}\,\int^t_0
\|\nabla\vphi\|^{1/\theta}_{L^\infty}\,\|\psi\|_{B^{s_1}_{p,r}}\,d\tau\,+ \nonumber \\
 & \quad+ \frac{C\,\eta}{\veps^{(1-\eta)/\eta}}\int^t_0
\|\nabla\psi\|^{1/\eta}_{L^\infty}\,\|\nabla\vphi\|_{B^{\sigma_1-1}_{p,r}}\,d\tau +
(1-\theta) \veps \|\psi\|_{\wtilde{L}^1_t(B^{s_2}_{p,r})} +
(1-\eta)\,\veps\,\|\nabla\vphi\|_{\wtilde{L}^1_t(B^{\sigma_2-1}_{p,r})}\,.
\label{conest:comm-deri}
\end{align}
\end{itemize}
\end{lemma}
The proof of \eqref{linest:comm} is classical while the proof of \eqref{conest:comm-deri} can be found in the appendix.
Let's just give a remark which will be used for estimating the density:
\begin{rem}\label{rem:comm}
If $s\geq d/p$ with $r=1$ when $s=d/p$, then  \eqref{conest:comm-deri} with $\theta=1/2$, $\eta=1$ becomes
\begin{equation*}\label{linest:comm,kappa}
 \Bigl\|2^{js}\int^t_0
\left\|\nabla\left(\left[\vphi,\Delta_j\right]\nabla\psi\right)\right\|_{L^p}\,d\tau
\Bigr\|_{\ell^r}
\,\leq
\,C_\eps\,\int^t_0
  \|\nabla\vphi\|_{B^s_{p,r}}^2\|\psi\|_{B^s_{p,r}}
   \,+\,
 \,\veps\,\|\psi\|_{\wtilde{L}^1_t(B^{s+2}_{p,r})} .
\end{equation*}
Indeed, one only has to use the interpolation inequality of the embedding $B^{d/p+2}_{p,1}\hookrightarrow C^{0,1}\hookrightarrow B^{d/p}_{p,1}$ (noticing also $B^{d/p}_{p,1}\hookleftarrow B^s_{p,r}$, $L^1_t(B^{d/p+2}_{p,1})\hookleftarrow \tilde L^1_t(B^{s+2}_{p,r})$) and Young's inequality.
\end{rem}

Next lemma is in the same spirit of \eqref{conest:comm-deri}, by view of Proposition \ref{p:prod} which gives estimates for  the product of two functions.
 The proof is quite similar (and  easier) and hence omitted.
\begin{lemma} \label{l:prod}
 Let $s>0$, $s_i, \sigma_i\in \R$, $i=1,2$, $(p,r)\in[1,+\infty]^2$, such that  $s=\theta s_1+(1-\theta)s_2=\eta\sigma_1+(1-\eta)\sigma_2$ with $\theta, \eta\in (0,1]$.
Then the following holds:
\begin{eqnarray*}
  \|f\,g\|_{\wtilde{L}^1_t(B^s_{p,r})} & \leq & \frac{C\,\theta}{\veps^{(1-\theta)/\theta}}\,\int^t_0
\|f\|^{1/\theta}_{L^\infty}\,\|g\|_{B^{s_1}_{p,r}}\,d\tau
\,+\,\frac{C\,\eta}{\veps^{(1-\eta)/\eta}}\int^t_0
\|g\|^{1/\eta}_{L^\infty}\,\|f\|_{B^{\sigma_1}_{p,r}}\,d\tau\,+ \\[1ex]
& & \qquad\qquad\qquad +\,(1-\theta)\,\veps\,\|f\|_{\wtilde{L}^1_t(B^{s_2}_{p,r})}\,+\,
(1-\eta)\,\veps\,\|g\|_{\wtilde{L}^1_t(B^{\sigma_2}_{p,r})}\,,
\forall \veps >0.
\end{eqnarray*}
\end{lemma}


Let us come  to  the linearized density equation
\begin{equation} \label{eq:l-vrho}
 \d_t\vrho\,+\,u\cdot\nabla\vrho\,-\,\div\left(\kappa\,\nabla\vrho\right)\,=\,f\,,
 \quad \vrho_{|t=0}=\vrho_0,
\end{equation}
for which one has the following a priori estimate:
\begin{prop} \label{p:lin_rho}
Let the triple $(s,p,r)$ verify
\begin{equation}\label{linearindex:rho}
s\geq \frac dp\,,\;\;p\,\in\,(1,+\infty)\,,\;\;r\in[1,+\infty]\,,\qquad\;\hbox{ with }\qquad\;
r=1\quad\textrm{ if }\quad s=1+\frac dp\;\textrm{ or }\;\frac dp.
\end{equation}
Let $u, \kappa, f$ be   smooth  such that $\nabla u\in B^{s-1}_{p,r}$,   $\kappa\geq \kappa_\ast>0$, $\nabla\kappa\in B^s_{p,r}$ and $f\in\wtilde{L}^1_{T_0}(B^s_{p,r})$.
{}
Then there exists a positive constant $C_1$ (depending only on $\kappa_\ast, d,s,p,r$) such that, for every smooth solution $\vrho$ of
\eqref{eq:l-vrho}, the following estimate holds true
for every $t\in[0,T_0]$:
\begin{equation*} \label{est:l-vrho}
 \left\|\vrho\right\|_{\wtilde{L}^\infty_t(B^s_{p,r})}\,+\,\left\|\vrho\right\|_{\wtilde{L}^1_t(B^{s+2}_{p,r})}\;\leq\;
C_1\,e^{C_1 K(t)}\left(\|\vrho_0\|_{B^s_{p,r}}
\,+\,\|f\|_{\wtilde{L}^1_t(B^s_{p,r})}\right),
\end{equation*}
where we have defined $K(0)=0$ and
\begin{equation*}\label{linearest:K}
K'(t)\,:=\,1
+\|\nabla u\|_{B^{d/p}_{p,1}\cap B^{s-1}_{p,r}}
   \,+\,\|\nabla\kappa \|_{B^{s}_{p,r}}^2\,.
\end{equation*}
\end{prop}


The proof is quite standard (see eg. the proof of Proposition 4.1 in \cite{D-L}): we apply the operator $\dj$ to the equation, we integrate first in space
and then in time; then we use the commutator estimates and Gronwall's Inequality to get the result.
Let's just sketch the proof.

Applying   $\Delta_j$ to Equation \eqref{eq:l-vrho} yields
\begin{equation} \label{eq:th-loc}
\d_t\vrho_j\,+\,u\cdot\nabla\vrho_j\,-\,\div\left(\kappa\,\nabla\vrho_j\right)\,=\,f_j\,+\,\mc R^1_j\,-\,\mc R^2_j\,,
\end{equation}
where we have set $\vrho_j\,:=\,\Delta_j\vrho$, $f_j\,:=\,\Delta_jf$, $\mc R^1_j\,:=\,\left[u,\Delta_j\right]\cdot\nabla\vrho$ and
$\mc R^2_j\,:=\,\div\left[\kappa,\Delta_j\right]\,\nabla\vrho$.

Hence one has
$$
\frac{d}{dt}\|\vrho_j\|^p_{L^p}\,+\,C\,2^{2j}\,\|\vrho_j\|^p_{L^p}\,\leq\,C\,\|\vrho_j\|^{p-1}_{L^p}
\left(\left\|f_j\right\|_{L^p}\,+\,\left\|\mc R^1_j\right\|_{L^p}\,+\,\left\|\mc R^2_j\right\|_{L^p}\right)\,,\quad j\geq 0,
$$
thanks to  the following Bernstein type inequality (see e.g.  Appendix B of \cite{D})
\begin{equation*}\label{ineq:Bernstein}
-\Int_{\R^d}  \div(\kappa\nabla \vrho_j) |\vrho_j|^{p-2} \vrho_j
 =(p-1) \Int_{\R^d} \kappa |\nabla\varrho_j|^2 |\varrho_j|^{p-2}
 \geq C(d,p,\kappa_\ast)\, 2^{2j} \Int_{\R^d} |\varrho_j|^p,\,\quad j\geq 0,\,p\in (1,\infty).
\end{equation*}
Since the equality in the above  holds also for $j=-1$,  the following holds:
\begin{align}
& \left\|\vrho \right\|_{\tilde L^\infty_t(B^s_{p,r})}\,+\,\left\|\vrho\right\|_{\wtilde{L}^1_t(B^{s+2}_{p,r})}
  \leq\,
C\biggl(\|\vrho_0\|_{B^s_{p,r}}
       \,+\,2^{-(s+2)}\,\left\|\Delta_{-1}\vrho\right\|_{L^1_t(L^p)}\,+
       \nonumber\\ 
& \qquad\qquad\qquad
        \,+\,\|f\|_{\wtilde{L}^1_t(B^s_{p,r})}
        \,+\,  \Bigl\|2^{js}\Int^t_0\left\|\mc R^1_j(\tau)\right\|_{L^p}d\tau\Bigr\|_{\ell^r}
       \,+\,\Bigl\|2^{js}\Int^t_0\left\|\mc R^2_j(\tau)\right\|_{L^p}d\tau\Bigr\|_{\ell^r}
       \biggr)\,.\label{est:vrho_part}
\end{align}
The low-frequency term $\Delta_{-1}\vrho$ can be easily bounded in $[0,T_0]$:
\begin{equation*} \label{est:low-f}
2^{-(s+2)}\,\left\|\Delta_{-1}\vrho\right\|_{L^1_t(L^p)}\;\leq\;C\,\int^t_0\|\vrho\|_{B^{s}_{p,r}}\,d\tau\,.
\end{equation*}
One applies Proposition \ref{p:comm} on   the first commutator term and Remark \ref{rem:comm}  on the second commutator term.
Finally performing Gronwall's inequality on \eqref{est:vrho_part} gives the conclusion.

\begin{rem}
Let us point out here that in the process of proving the uniqueness of the solutions to System \eqref{system}, there is one derivative loss for the difference of two solutions.
 We therefore  have to look for a priori estimates for the unknowns in $B^s_{p,r}$, under a weaker condition \eqref{linearindex:rho} on the indices (instead of \eqref{index:s,r}).
\end{rem}


The linearized equation for the velocity reads
\begin{equation}\label{lineq:u}
\left\{
\begin{array}{cc}
&\d_t u+w\cdot\nabla u+\lambda\nabla\pi=h,\\
&\div u=0,\\
&u|_{t=0}=u_0,
\end{array}
\right.
\end{equation}
where the initial datum $u_0$, the transport vector field $w$, the coefficient $\lambda$ and the source term $h$ are all smooth
and decrease rapidly at infinity, such that $ \lambda\geq \lambda_\ast>0$. We have the following a priori estimate:
\begin{prop}\label{p:TDnu}
Let
\begin{equation}\label{linearindex:u}
s>\dfrac dp-\dfrac d4,\quad p\in [2,4],\quad r\in[1,+\infty]\,,\qquad\hbox{ with }\qquad
r=1\;\textrm{ if }\;s=\frac dp\;\textrm{ or }\;1+\frac dp\,.
\end{equation}
Then   the following estimates
 hold true:
\begin{align}
&\|u\|_{\tilde L^\infty_t(B^{s}_{p,r})}  \leq
C_2e^{C_2 W(t)}\Bigl(\|u_0\|_{B^{s}_{p,r}}+\|h\|_{\tilde L^1_t(B^{s}_{p,r})\cap L^1_t(L^2)}\Bigr),
\quad  W(t)=\int^t_0 \|\nabla w \|_{B^{\frac dp}_{p,1}\cap B^{s-1}_{p,r}},\label{linearest:u} \\
&\|\nabla\pi\|_{\tilde L^1_{t}(B^s_{p,r})\cap L^1_{t}(L^2)}
  \leq   C_2\Bigl(\|h\|_{\tilde L^1_t(B^{s}_{p,r})\cap L^1_t(L^2)}
+W(t)\|u\|_{\tilde L^\infty_t(B^s_{p,r})}\Bigr),\label{linearest:pi}
\end{align}
where      $C_2=C_2(t)$   is  a positive \textit{time-dependent function}, depending only on $d,p,s,r,\lambda_\ast,\lambda^\ast(T)$, with
$$\lambda^\ast(t):=\|\lambda\|_{L^\infty_t(L^\infty)}+
 \|\nabla\lambda\|_{\tilde L^\infty_t(B^{\frac dp}_{p,1}\cap B^{s-1}_{p,r})}.$$
\end{prop}
The proof is similar as in \cite{D} and let us also just sketch it.\\
Firstly, Proposition \ref{p:comm} entails  the following estimate for $u$ :
\begin{align}\label{linest:u}
\|u(t)\|_{\tilde L^\infty_t(B^{s}_{p,r})}\leq
 \Bigl(\|u_0\|_{B^{s}_{p,r}}+\Int^t_0 W' \| u\|_{B^s_{p,r}}+\|h-\lambda\nabla \pi\|_{\tilde L^1_t(B^{s}_{p,r})}\Bigr).
\end{align}
By product estimates in Proposition \ref{p:prod}, we have (noticing  $\lambda\nabla\pi\equiv
(\Delta_{-1}\lambda)\nabla\pi+ ((Id-\Delta_{-1})\lambda  ) \nabla \pi$)
\begin{align*} \label{est:l-Pi}
\|\lambda\nabla\pi\|_{\tilde L^1_t(B^s_{p,r})}\leq  C \lambda^\ast\|\nabla\pi\|_{\tilde L^1_t(B^s_{p,r})},
\hbox{ if }s>-\min\{\frac{d}{p},\frac{d}{p'}\}, \hbox{ with }r=1 \hbox{ if }s=\frac{d}{p}.
\end{align*}


Thanks to  $\div u\equiv 0$, $ \pi $ satisfies the following elliptic equation:
\begin{equation*}\label{lineareq:pi,origin}
\div(\lambda\nabla\pi)=\div(h-w\cdot\nabla u)=\div(h-u\cdot\nabla w+u\,\div w).
\end{equation*}
 Similar as to get \eqref{est:vrho_part}, one finds
\begin{align}\label{linest:pi,origin}
 \|\nabla\pi\|_{\tilde L^1_t(B^s_{p,r})}
&\lesssim \|\nabla\Delta_{-1}\pi\|_{L^1_t(L^p)}+\|\div (h-w\cdot\nabla u)\|_{\tilde L^1_t(B^{s-1}_{p,r})} +\|2^{j(s-1)}\|\div[\lambda,\dj]\nabla\pi\|_{L^1_t(L^p)}\|_{\ell^r}.
\end{align}
To bound the above  commutator term, one applies  Proposition \ref{p:comm} to it;
 then one uses the following  interpolation
inequality (with some $\epsilon, \eta \in (0,1)$)
\begin{align*}
\|\nabla\pi\|_{L^1_t(B^{s-1}_{p,r})}
&\lesssim \|\nabla\pi\|_{\tilde L^1_t(B^{s+\epsilon-1}_{p,r})}
\lesssim\|\nabla\pi\|_{\tilde L^1_t(B^{\frac dp-\frac d2}_{p,\infty})}^{1-\eta}\|\nabla\pi\|_{\tilde L^1_t(B^s_{p,r})}^{\eta}
\lesssim \|\nabla\pi\|_{ L^1_t(L^2)}^{1-\eta}\|\nabla\pi\|_{\tilde L^1_t(B^s_{p,r})}^{\eta};
\end{align*}
  finally  one has
\begin{align*}
\|2^{j(s-1)}\|\div[\lambda,\dj]\nabla\pi\|_{L^1_t(L^p)}\|_{\ell^r}
\leq C(d,s,p,r,\eps,\epsilon,\lambda^\ast)\|\nabla\pi\|_{L^1_t(L^2)} +\eps \|\nabla\pi\|_{\tilde L^1_t(B^s_{p,r})}.
\end{align*}
Therefore, thanks to $\|\nabla\Delta_{-1}\pi\|_{L^p}\lesssim \|\nabla\pi\|_{L^2}$ and Lemma \ref{l:laxmilgram},    it rests to dealing with
$$
\|u\cdot\nabla w\|_{L^1_t(L^2)},\quad \|u\,\div w\|_{L^1_t(L^2)} \quad \hbox{and}\quad \|\div (w\cdot\nabla u)\|_{\tilde L^1_t(B^{s-1}_{p,r})}.
$$
For $p\leq 4$, $s> d/p- d/4$, we can easily find that
\begin{equation*}\label{linest:w,u,L2}
\|u\,\div w\|_{L^1_t(L^2)},\,\|u\cdot\nabla w\|_{L^1_t(L^2)}\leq \Int^t_0 \|u\|_{L^4} \|\nabla w\|_{L^4}
\lesssim
\Int^t_0 \|u\|_{B^{s}_{p,r}}\|\nabla w\|_{B^{\frac dp}_{p,\infty}}\,d\tau\,.
\end{equation*}
On the other hand,  it is easy to decompose
$\|\div(w\cdot\nabla u)\|_{B^{s-1}_{p,r}}$ into
\begin{align*}
\|T_{\d_i w^j}\d_j u^i+T_{\d_j u^i}\d_i w^j+
\div\bigl(R(w^j,\d_j u)\bigr)\|_{B^{s-1}_{p,r}},
\end{align*}
which can be controlled, according to Proposition \ref{p:prod}, by
$  W'(t)\|\nabla u\|_{B^{s-1}_{p,r}}$.


To conclude, Estimate \eqref{linearest:pi} holds, and so does estimation \eqref{linearest:u}, by view of \eqref{linest:u} and Gronwall's Inequality.

%

\subsection{Proof of the existence}\label{subs:proof_th1}

In this subsection we will follow the standard procedure to prove the local existence of the solution to  System
\eqref{system}: we construct a sequence of approximate solutions which have uniform bounds and then we prove the convergence to a unique solution.
{}
In particular, in order to bound the nonlinearities,  the density should be \textit{small} when \textit{integrated in time}.
Since we admit also large initial density $\rho_0$, we will introduce the large linear part $\rho_{L} $ of
the solution $\rho$, so that the remainder part $\bar\rho:=\rho-\rho_L$ is small and hence easier to handle.
{}
 In the convergence part,  we will first show convergence  in a space with \textit{lower regularity} (i.e. in space $E^{d/p}_{p,r}(T)$, see \eqref{space:E}) and
then the solution is in $E^s_{p,r}$ by Fatou's property.

 We will freely use the following estimates (by Propositions \ref{p:prod} and \ref{p:comp_grad}):
     \begin{equation}\label{ineq}
  \|uv\|_{B^s_{p,r}}\lesssim \|u\|_{B^s_{p,r}}\|v\|_{B^s_{p,r}},\,
\|f(\rho) \|_{B^s_{p,r}}\leq C(\|\rho\|_{L^\infty})\|\varrho\|_{B^s_{p,r}}
 \hbox{ with }f(1)=0,\,s>\frac{d}{p}\hbox{ or }s\geq \frac{d}{p}, r=1,
 \end{equation}
and their time-dependent version
\begin{align*}
 & \|uv\|_{\tilde L^q_t(B^s_{p,r})}
  \lesssim \|u\|_{\tilde L^{q_1}_t(B^s_{p,r})}
          \|v\|_{\tilde L^{q_2}(B^s_{p,r})},\,
 \|f(\rho) \|_{\tilde L^q_t(B^s_{p,r})}
\leq C( \|\rho\|_{L^\infty_t(L^\infty)})\|\varrho\|_{\tilde L^q_t(B^s_{p,r})}
\quad
          \hbox{ with }\frac 1q=\frac{1}{q_1}+\frac{1}{q_2} .
\end{align*}


\subsubsection{Step 1 -- Construction of a sequence of approximate solutions}

In this step, we take $(s,p,r)$ such that Conditions \eqref{index:s,r} and \eqref{index:p} hold true. Let us introduce the approximate
solution sequence $ \{(\varrho^{n}, u^{n}, \nabla\pi^{n})\}_{n\geq 0} $ by induction.

Without loss of generality we can assume
\begin{equation*}\label{seqest:initial}
\frac{\rho_\ast}{2}\,\leq\, S_{n}\rho_0,\qquad \forall n\in \N\,;
\end{equation*}
then, first of all we set $(\varrho^{0}, u^{0}, \nabla\pi^{0}):=(S_{0}\varrho_{0},S_{0}u_{0},0)$, which are  smooth and fast decaying at infinity.

Now, we assume by induction that the triplet $(\varrho^{n-1}, u^{n-1}, \nabla \pi^{n-1})$
of smooth and fast decaying functions has been constructed.
Besides, let us suppose also that there exists a sufficiently small
parameter $\tau$ (to be determined later),
a positive time $T^{*}$ (which may depend on $\tau$)  and
a positive constant $C_M$ (which may depend on $M$) such that
\begin{align}
&\frac{\rho_\ast}{2}\leq \rho^{n-1}:=1+\varrho^{n-1}\,,
\quad \|\varrho^{n-1}\|_{\tilde L^{\infty}_{T^{*}}(B^{s}_{p,r})}\leq C_M\,,
\quad
\|\varrho^{n-1}\|_{ \tilde L^{2}_{T^{*}}(B^{s+1}_{p,r})\cap \tilde L^{1}_{T^{*}}(B^{s+2}_{p,r})}
       \leq \tau\,,\label{seqest:rho}\\[1ex]
&\|u^{n-1}\|_{\tilde L^{\infty}_{T^{*}}(B^{s}_{p,r})}\leq C_M,\;
\|u^{n-1}\|_{L^{2}_{T^{*}}(B^{s}_{p,r})\cap L^{1}_{T^{*}}(B^{s}_{p,r})}\leq \tau,
\;\|\nabla\pi^{n-1}\|_{\tilde L^1_{T^\ast}(B^s_{p,r})\cap L^1_{T^\ast}(L^2)}\leq \tau^{1/2}.\label{seqest:u}
\end{align}
Remark that the above estimates \eqref{seqest:rho} and \eqref{seqest:u} obviously
hold true for $(\varrho^{0}, u^{0},\nabla\pi^{0})$, if $T^{*}$ is assumed to be small enough.

Now we define $ (\varrho^{n}, u^{n}, \nabla\pi^{n})$   as the unique smooth global solution  of the linear system
\begin{equation}\label{seqeq}
\left \{\begin{array}{cc}
&\d_{t}\varrho^{n}+u^{n-1}\cdot\nabla\varrho^{n}-\div(\kappa^{n-1}\nabla\varrho^{n})=0,\\
&\d_{t}u^{n}+(u^{n-1}+\nabla b^{n-1})\cdot\nabla u^{n}+\lambda^{n-1}\nabla\pi^{n}=h^{n-1},\\
&\div u^{n}=0,\\
&(\varrho^{n}, u^{n})|_{t=0}=(S_{n}\varrho_{0},S_{n}u_{0}),
\end{array}
\right .
\end{equation}
where we have set $a^{n-1}=a(\rho^{n-1})$,
$b^{n-1}=b(\rho^{n-1})$, $\kappa^{n-1}=\kappa(\rho^{n-1})$, $\lambda^{n-1}=\lambda(\rho^{n-1})$  and
\begin{align}\label{L2seqeq:h}
h^{n-1}&= (\rho^{n-1})^{-1}\div((u^{n-1}+\nabla b^{n-1})\otimes\nabla a^{n-1})
 \nonumber \\
 &=(\rho^{n-1})^{-1}\Bigl( \Delta b^{n-1} \nabla a^{n-1}
    \,+\,  u^{n-1}\cdot \nabla^2 a^{n-1}
   \, + \,\nabla b^{n-1}\cdot\nabla^2 a^{n-1} \Bigr) .
\end{align}
We want to show that also the triplet $ (\varrho^{n}, u^{n}, \nabla\pi^{n})$ verifies \eqref{seqest:rho} and \eqref{seqest:u}.

First of all,  we apply the maximum principle to the linear parabolic equation for
$\varrho^n$, yielding $\rho^n:=1+\vrho^n\in [\rho_\ast/2, \rho^\ast]$.
{}
Now, we  introduce  $\varrho_{L}$ as the solution of the heat equation with the initial datum $ \varrho_{0}\in B^s_{p,r} $:
$$
\left\{\begin{array}{l}
        \d_t\vrho_L\,-\,\Delta\vrho_L\,=\,0 \\[1ex]
	\left(\vrho_L\right)_{|t=0}\,=\,\vrho_0\,.
       \end{array}\right.
$$
Then,  for any positive time $T<+\infty$, there exists   some constant $C_T>0$ depending on $T$ such that
\begin{equation}\label{seqest:rho,L}
\|\varrho_{L}\|_{\tilde L^{\infty}_{T }(B^{s}_{p,r})}\,+\,\|\varrho_{L}\|_{\tilde L^{1}_{T }(B^{s+2}_{p,r})}\;\leq\;
C_T\,\|\varrho_0\|_{B^s_{p,r}}\,.
\end{equation}
Furthermore, given $\tau>0$, we can choose $T^\ast<+\infty$ such that one has
 \begin{equation}\label{seqest:rho,L,tau}
 \|\varrho_L\|_{\tilde  L^{2}_{T^\ast}(B^{s+1}_{p,r})\cap\tilde L^{1}_{T^\ast}(B^{s+2}_{p,r})}\leq \tau^2.
 \end{equation}
 Indeed, by definition we have
 \begin{align*}
  \|\varrho_L\|_{\tilde L^{1}_{T^\ast}(B^{s+2}_{p,r})}
  &=\Bigl\|\Bigl(2^{js }\Int^{T^*}_0 \|2^{2j}e^{t\Delta}\Delta_j\varrho_0\|_{L^p}\,dt\Bigr)_j\Bigr\|_{\ell^r}.
 \end{align*}
The operator $e^{t\Delta}\Delta_j$ belongs to $\mc{L}(L^p):=\left\{A:L^p\ra L^p\;\mbox{linear and bounded}\,\right\}$: more precisely,
$$
\left\|e^{t\Delta}\,\Delta_j\right\|_{\mc{L}(L^p)}\,\leq\,C \;\forall j\geq -1,
\quad\hbox{and}\quad
\left\|e^{t\Delta}\,\Delta_j\right\|_{\mc{L}(L^p)}\,\leq\,C\,e^{-C\,t\,2^{2j} }\;\forall j\geq 0.
$$
Then, for some fixed $N$ large enough, we infer
$$
\Bigl\| \Bigl(\Int^{T^\ast}_0 2^{2j} e^{-C\,t\,2^{2j} } \,dt\Bigr)_{0\leq j\leq N} \Bigr\|_{\ell^\infty}\leq C(1-e^{-C 2^{2N} T^\ast });
$$
From this,  by decomposing   $\vrho_0$ into low frequencies  (large part)   and high frequencies (small part) and choosing
$T^\ast$ small enough, one gathers $\|\varrho_L\|_{\tilde L^{1}_{T^\ast}(B^{s+2}_{p,r})}\leq \tau^2$.
The term $\|\vrho_L\|_{\tilde  L^2_{T^*}(B^{s+1}_{p,r})}$ can be handled in the same way or by interpolation inequality.
Hence, our claim  \eqref{seqest:rho,L,tau} is proved.

Now we define the sequence $ \varrho_{L}^{n}=S_{n}\varrho_L $: it too solves the free heat equation,
but with initial data $S_{n}\varrho_{0}$. Hence, it too satisfies \eqref{seqest:rho,L} and \eqref{seqest:rho,L,tau}.


 We next consider the small remainder $ \bar{\varrho}^{n}\,:=\,\varrho^{n}-\varrho^{n}_{L} $. We claim that it fulfills,
for all $n\in\N$,
\begin{equation}\label{seqest:theta,bar}
\|\bar{ \varrho }^{n}\|_{\tilde L^{2}_{T^{*}}(B^{s+1}_{p,r})}\leq \|\bar{ \varrho }^{n}\|_{\tilde L^{ \infty }_{T^{*}}(B^{s}_{p,r})}+\Vert \bar{\varrho}^{n}  \Vert_{\tilde L^{1}_{T^{*}}(B^{s+2}_{p,r})}\leq \tau^{3/2}\,.
\end{equation}
In fact, $ \bar{\varrho}^{n}=\varrho^{n}-\varrho^{n}_{L} $ solves
\begin{equation}\label{seqeq:theta,bar}
\left \{
\begin{array}{cc}
&\d_{t}\bar{\varrho}^{n}+u^{n-1}\cdot\nabla \bar{\varrho}^{n}-
\div (\kappa^{n-1} \nabla \bar{\varrho}^{n})=-u^{n-1}\cdot\nabla \varrho^{n}_{L}+\div ((\kappa^{n-1}-1)\nabla\varrho^{n}_{L}),\\
&\bar{\varrho}^{n}|_{t=0}=0.
\end{array}
\right .
\end{equation}
So, if we define
$$
K^{n-1}(t)\;:=\;t+\left\|\nabla u^{n-1}\right\|_{L^1_t(B^{s-1}_{p,r})}
\,+\,
\left\|\nabla\kappa^{n-1}\right\|_{L^2_t(B^{s}_{p,r})}^2 ,
$$
by Proposition \ref{est:l-vrho} we infer that
$$
\|\bar\varrho^n\|_{\tilde L^\infty_{T^{*}}(B^s_{p,r})\cap \tilde L^1_{T^\ast}(B^{s+2}_{p,r})}\,\leq\,C\,e^{C K^{n-1}(T^*)}\,
\left\|-u^{n-1}\cdot\nabla \varrho^{n}_{L}+\div \Bigl((\kappa^{n-1}-1)\nabla\varrho^{n}_{L}\Bigr)\right\|_{\wtilde{L}^1_{T^*}(B^s_{p,r})}\,.
$$
Inductive assumptions and estimate \eqref{seqest:rho,L,tau} for $\varrho_L$ help to bound the right-hand side in the above relation
by  $  CC_M\tau^2$. Therefore, \eqref{seqest:theta,bar} is proved, and hence \eqref{seqest:rho} holds for
$\vrho^n\,=\,\oline{\vrho}^n+\vrho^n_L$, for sufficiently small $\tau$.


We now want  to get \eqref{seqest:u}, relying mainly on  Proposition \ref{p:TDnu}.
In fact, product estimates \eqref{ineq} and the embedding result in Remark \ref{rem:besov} entail (noticing also $L^4\hookleftarrow B^{s-1}_{p,\infty}$)
$$
 \Vert h^{n-1}\Vert_{\tilde L^{1}_{T^\ast}(B^{s}_{p,r})\cap L^1_{T^\ast}(L^2)},\quad
W^{n-1}(T^\ast):=\Int^{T^\ast}_0\|\nabla u^{n-1}+\nabla^2 b^{n-1}\|_{B^{s-1}_{p,r}}
\leq C\tau.
$$
Thus, applying Proposition \ref{p:TDnu} to system \eqref{seqeq} implies
\begin{align*}
 \Vert u^{n} \Vert_{\tilde L^{\infty}_{T^\ast}(B^{s}_{p,r})}
\leq C(T^\ast)(\Vert  S_{n}u_{0}\Vert_{B^{s}_{p,r}}
+C\tau)\leq C_M,\quad
\|\nabla\pi^n\|_{\tilde L^1_{T^\ast}(B^s_{p,r})\cap L^1_{T^\ast}(L^2)}
\leq C\tau+C_MC\tau.
\end{align*}
Hence   \eqref{seqest:u} also holds true for small $\tau$ and $T^\ast$.

\subsubsection{Step 2 -- Convergence of the sequence} \label{subs:step2}

In this step we will consider the ``difference'' sequence
$$
(\delta \varrho^{n},\delta u^{n},\nabla\delta\pi^{n}) := (\varrho^{n}-\varrho^{n-1},u^{n}-u^{n-1},\nabla\pi^{n}-\nabla\pi^{n-1}),
\qquad\forall\,n\geq1\,
$$
in the Banach space $E^{d/p}_{p,1}(T^\ast)$ (recall \eqref{space:E} for its definition).

First of all, by System \eqref{seqeq}, $(\delta \varrho^{n},\delta u^{n},\nabla\delta\pi^{n})$ solves
\begin{equation}\label{seqeq:diff}
\left \{
\begin{array}{cc}
&\d_{t}\delta \varrho^{n}+u^{n-1}\cdot\nabla\delta \varrho^{n}-\div(\kappa^{n-1}\nabla\delta \varrho^{n})=F^{n-1},\\
&\d_{t}\delta u^{n}+(u^{n-1}+\nabla b^{n-1})\cdot\nabla\delta u^{n}+\lambda^{n-1}\nabla\delta\pi^{n}=H^{n-1},\\
&\div\delta u^{n}=0,\\
&(\delta \varrho^{n},\delta u^{n})|_{t=0}=(\Delta_{n}\varrho_{0},\Delta_{n}u_{0}),
\end{array}\right .
\end{equation}
where
\begin{align*}
& F^{n-1}=-\delta u^{n-1}\cdot\nabla\varrho^{n-1}+\div((\kappa^{n-1}-\kappa^{n-2})\nabla\varrho^{n-1}),\\
& H^{n-1}=h^{n-1}-h^{n-2}-(\delta u^{n-1}+\nabla\delta b^{n-1})\cdot\nabla u^{n-1}-(\lambda^{n-1}-\lambda^{n-2})\nabla\pi^{n-1}.
\end{align*}

Next we apply a priori estimates given by Propositions \ref{p:lin_rho} and \ref{p:TDnu},  with $ s=d/p$, $p\in [2,4]$ and
$r=1 $, to $\delta\varrho^{n}$ and $ (\delta u^{n},\nabla\delta\pi^{n}) $ respectively.
 The use of inductive assumptions   gives us
\begin{align}
&\|\delta\varrho^{n}\|_{L^{\infty}_{T^\ast}(B^{d/p}_{p,1})\cap L^{1}_{T^\ast}(B^{d/p+2}_{p,1})}\,\leq\,
C\left (\|\Delta_{n}\varrho_{0}\|_{B^{d/p}_{p,1}}+\|F^{n-1}\|_{L^{1}_{T^\ast}(B^{d/p}_{p,1})}\right ),\label{varseqest:rho}\\
&\|\delta u^{n}\|_{\tilde L^{\infty}_{T^\ast}(B^{d/p}_{p,1})}
+\|\nabla\delta\pi^{n}\|_{ L^{1}_{T^\ast}(B^{d/p}_{p,1}\cap L^2)}
\leq C \Bigl(\|\Delta_{n}u_{0}\|_{B^{d/p}_{p,1}}+\|H^{n-1}\|_{L^{1}_{T^\ast}(B^{d/p}_{p,1}\cap L^2)} \Bigr ).
\label{varseqest:u}
\end{align}
Next we use the following fact, coming from Proposition \ref{p:comp_grad}: for $f$ smooth, $s\geq d/p$,
$$
\|\delta f^{m}\|_{B^{s}_{p,1}}
:=\|f(\rho^m)-f(\rho^{m-1})\|_{B^{s}_{p,1}}
\,\leq\,C(\|\varrho^{m}\|_{B^{s}_{p,1}},\|\varrho^{m-1}\|_{B^{s}_{p,1}})\;
\|\delta\varrho^{m}\|_{B^{s}_{p,1}}.
$$
Therefore, one easily gets
\begin{align*}
&\|F^{n-1}\|_{L^1_{T^\ast}(B^{d/p}_{p,1})}
\leq C (\|\delta\varrho^{n-1}\|_{L^2_{T^\ast}(B^{d/p+1}_{p,1})}
\|\varrho^{n-1}\|_{L^2_{T^\ast}(B^{d/p+1}_{p,1})}\\
&\qquad\qquad\qquad\qquad +\Int^{T^\ast}_0 \|\delta u^{n-1}\|_{B^{d/p}_{p,1}}\|\nabla\varrho^{n-1}\|_{B^{d/p}_{p,1}}
+\|\delta\varrho^{n-1}\|_{B^{d/p}_{p,1}}\|\varrho^{n-1}\|_{B^{d/p+2}_{p,1}}),\\
&\|H^{n-1}\|_{L^1_{T^\ast}(B^{d/p}_{p,1}\cap L^2)}
\leq C (\|\delta\varrho^{n-1}\|_{L^1_{T^\ast}(B^{d/p+2}_{p,1})}\\
&\qquad\qquad+\|\delta\varrho^{n-1}\|_{L^2_{T^\ast}(B^{d/p+1}_{p,1})}
\Bigl(\|(\varrho^{n-1},\varrho^{n-2})\|_{L^2_{T^\ast}(B^{d/p+2}_{p,1})}
+\|(u^{n-1},u^{n-2})\|_{L^2_{T^\ast}(B^{d/p+1}_{p,1})}\Bigr)\\
&\qquad\qquad
+\Int^{T^\ast}_0 \|\delta\varrho^{n-1}\|_{B^{d/p}_{p,1}}\|\nabla\pi^{n-1}\|_{B^{d/p}_{p,1}}
+\|\delta u^{n-1}\|_{B^{d/p}_{p,1}}(\|\varrho^{n-1}\|_{B^{d/p+2}_{p,1}}+\|u^{n-1}\|_{B^{d/p+1}_{p,1}})).
\end{align*}
Plugging the uniform estimate \eqref{seqest:rho} into \eqref{varseqest:rho} to get bounds
on $\|\delta\varrho^{n-1}\|_{L^1_{T^\ast}(B^{d/p+2}_{p,1})}$ (which appears in $H^{n-1}$) and then relation \eqref{varseqest:u} becomes
\begin{align*}
&\|\delta u^{n}\|_{\tilde L^{\infty}_{T^\ast}(B^{d/p}_{p,1})}+\|\nabla\delta\pi^{n}\|_{\tilde L^{1}_{T^\ast}(B^{d/p}_{p,1})}\leq C (\|(\Delta_{n}u_{0},\Delta_{n-1}\varrho_{0})\|_{B^{d/p}_{p,1}}
+\tau\|(\delta\varrho^{n-1},\delta\varrho^{n-2})\|_{L^2_{T^\ast}(B^{d/p+1}_{p,1})}\\
&\quad\qquad\qquad
+\Int^{T^\ast}_0 \|(\delta \varrho^{n-1},\delta u^{n-1},\delta u^{n-2})\|_{B^{d/p}_{p,1}}
\Bigl(\|\nabla\pi^{n-1}\|_{B^{d/p}_{p,1}}
+\|(\varrho^{n-1},\varrho^{n-2})\|_{B^{d/p+2}_{p,1}}
+\|u^{n-1}\|_{B^{d/p+1}_{p,1}}\Bigr).
\end{align*}

Let us now define
$$
B^{n}(t)\,:=\,\|\delta\varrho^{n}\|_{L^{\infty}_{t}(B^{d/p}_{p,1})}\,+\,
\|\delta\varrho^{n}\|_{L^{1}_{t}(B^{d/p+2}_{p,1})}\,+\,\|\delta u^{n}\|_{L^{\infty}_{t}(B^{d/p}_{p,1})}\,+\,
\|\nabla\delta\pi^{n}\|_{L^{1}_{t}(B^{d/p}_{p,1}\cap L^2)}\,;
$$
then, from previous inequalities we gather
\begin{eqnarray*}
B^{n}(t)  \leq  C\left\|(\Delta_{n-1}\varrho_{0},
\Delta_{n}\varrho_{0},\Delta_{n}u_{0})\right\|_{B^{d/p}_{p,1}}+
\tau^{\frac 12}\left(B^{n-1}(t)+B^{n-2}(t)\right)
+C\Int^{t}_0\left(B^{n-1}+B^{n-2}\right)D(\sigma)d\sigma,
\end{eqnarray*}
with $\|D(t)\|_{L^1([0,T^\ast])}\leq C$.
Noticing
$$
\left\|\left(\Delta_n\vrho_0,\Delta_nu_0\right)\right\|_{B^{d/p}_{p,1}}\,\leq\,C\,2^{n(d/p)}\,
\left\|\left(\Delta_n\vrho_0,\Delta_nu_0\right)\right\|_{L^p}\,\forall n\geq 0,
$$
it follows $ \sum_n B^{n}(t)\,<\,+\infty $ uniformly in $ [0,T^{*}] $.
Hence, we gather that the sequence $ {(\varrho^{n},u^{n},\nabla\pi^{n})} $ is a Cauchy sequence in the functional space
$E^{d/p}_{p,1}(T^*)$.
Then, it converges  to some $(\vrho,u,\nabla\pi)$, which actually belongs
to the space $E^s_{p,r}(T^*)$ by Fatou property.
Hence, by interpolation, the convergence holds true in any intermediate space between $E^s_{p,r}(T^*)$ and $E^{d/p}_{p,1}(T^*)$,
and this is enough to pass to the limit in our equations.
Thus, $(\vrho,u,\nabla\pi)$ is actually a solution of System \eqref{system}.
{}

The proof of uniqueness is exactly analogous to the above convergence proof, and hence omitted.

\section{Proof of Theorems \ref{th:cc_r} and \ref{th:lifespan}}\label{s:conti-life}
In this section we aim to get a continuation criterion and a lower bound of the lifespan for the local-in-time solutions given by Theorem \ref{th:w-p}.
It is only a matter of repeating a priori estimates established previously, but in an  ``accurate'' way (we use $L^\infty$-norm instead of   $B^{s-1}_{p,r}$-norm) for obtaining the continuation criterion, whereas in a ``rough'' way (we use \eqref{lifeprod}, \eqref{lifeineq:L2} below) for bounding the lifespan from below.

\subsection{Proof of the continuation criterion}\label{s:conti}
Theorem \ref{th:cc_r} actually issues easily from the following fundamental lemma.
\begin{lemma} \label{l:cc}
Let $s>0$, $p\in (1,+\infty)$ and $r\in [1,+\infty]$. Let $(\rho,u,\nabla\pi)$ be a solution of \eqref{system} over
$[0,T[\,\times\mbb{R}^d$ such that the hypotheses   in Theorem \ref{th:cc_r} hold true.
If $T$ is finite, then one gets
\begin{equation}\label{hyp_cc:finite}
\|\rho-1\|_{\wtilde{L}^\infty_T(B^s_{p,r})\cap\wtilde{L}^1_T(B^{s+2}_{p,r})}\,+\,
\|u\|_{\wtilde{L}^\infty_T(B^s_{p,r})}\,+\,\|\nabla\pi(t)\|_{\wtilde{L}^1_T(B^s_{p,r}) }\;<\;+\infty\,.
\end{equation}
\end{lemma}
In fact, to prove Theorem \ref{th:cc_r} from Lemma \ref{l:cc} is quite standard:  once   \eqref{hyp_cc:finite} and    Conditions \eqref{index:s,r}, \eqref{index:p} hold
true,  then there exists a positive time $t_0$ (thanks to Theorem   \ref{th:w-p}) such that, for any $\wtilde{T}<T$, System
\eqref{system} with initial data $\left(\rho(\wtilde{T}),u(\wtilde{T})\right)$ has a unique solution until the time $\tilde T+t_0$.
Thus, if we take, for instance, $\wtilde{T}=T-(t_0/2)$, then we get a solution until the time $T+(t_0/2)$,
which is, by uniqueness, the continuation of $(\rho,u,\nabla\pi)$. Theorem \ref{th:cc_r} then follows.

{}Therefore, we focus only on Lemma \ref{l:cc}: one uses $L^\infty$-norm (instead of Besov norm) to establish a priori estimates.

Let us consider the density term.
Our starting point is \eqref{eq:l-vrho},  with  $f=0$.
One argues as in proving Proposition \ref{p:lin_rho}, but controls commutators $\cR^1_j$ and $\cR^2_j$ (see \eqref{eq:th-loc} for definition)
by use of Commutator Estimates \eqref{linest:comm} and \eqref{conest:comm-deri} (with $\theta=\eta=1/2$) instead.
 More precisely, keeping in mind that $\kappa=\kappa(\rho)$, we arrive at
$$\displaylines{
\int^t_0\left\|2^{js}\,\left\|\mc{R}^1_j\right\|_{L^p}\right\|_{\ell^r}\,d\tau\,\leq\,\int^t_0
\left(\|\nabla u\|_{L^\infty}\,\|\vrho\|_{B^s_{p,r}}\,+\,\|\nabla\vrho\|_{L^\infty}\,\|u\|_{B^s_{p,r}}\right)d\tau\,,
\cr
\left\|2^{js}\,\Int^t_0\left\|\mc{R}^2_j\right\|_{L^p}\right\|_{\ell^r}\,d\tau\,
\leq\,\frac{C}{\veps}\int^t_0
\|\nabla\vrho\|^2_{L^\infty}\,\|\vrho\|_{B^s_{p,r}}\,d\tau\,
+\,\veps\,\|\vrho\|_{\wtilde{L}^1_t(B^{s+2}_{p,r})}\,.
}$$
Hence, \eqref{est:vrho_part} becomes
\begin{equation}\label{est_cc:vrho}
\|\vrho\|_{\wtilde{L}^{\infty}_t(B^s_{p,r})\cap\wtilde{L}^1_t(B^{s+2}_{p,r})}  \lesssim \left\|\vrho_0\right\|_{B^s_{p,r}}
 + \int^t_0 (1+\|\nabla u\|_{L^\infty} + \|\nabla\vrho\|^2_{L^\infty}) \|\vrho\|_{B^s_{p,r}} +
\int^t_0\|\nabla\vrho\|_{L^\infty} \|u\|_{B^s_{p,r}}   .
\end{equation}


Let us now consider velocity field and pressure term: we use Lemma \ref{l:tilde} to control the commutator
$\mc R_j\,:=\,\left[u+\nabla b(\rho)\,,\,\Delta_j\right]\cdot\nabla u$ and arrive at
\begin{align} \label{est:til_u_rough}
 \|u\|_{\wtilde{L}^\infty_t(B^s_{p,r})}\,\leq\,C\Bigl(\|u_0\|_{B^s_{p,r}}
+\int^t_0\|(\nabla u,\nabla^2\varrho)\|_{L^\infty}\,\|u\|_{B^s_{p,r}}
&+ \|\nabla u\|^2_{L^\infty}\,\|\vrho\|_{B^s_{p,r}}\,d\tau\,+\,
\eps \|\vrho\|_{\wtilde{L}^1_t(B^{s+2}_{p,r})}\nonumber\\
&+\|h\|_{\wtilde{L}^1_t(B^s_{p,r})}
+\left\|\lambda\,\nabla\pi\right\|_{\wtilde{L}^1_t(B^s_{p,r})}\Bigr).
\end{align}
Lemma \ref{l:prod} and Proposition \ref{p:comp_grad}   help us to control   the non-linear term  $h$ (see \eqref{h} for the definition):
\begin{align*}
 \|h\|_{\wtilde{L}^1_t(B^s_{p,r})}  \lesssim
 & \int^t_0
 \Bigl(\left\|\nabla^2\vrho\right\|_{L^\infty}+\|\nabla \varrho\|_{L^\infty}^2\Bigr)\,
\|u\|_{B^s_{p,r}}\,d\tau
\,+\,\Bigl(1+\|(\nabla\rho, u)\|_{L^\infty_t(L^\infty)}\Bigr)\,\|\vrho\|_{\wtilde{L}^1_t(B^{s+2}_{p,r})} \\
& + \int^t_0
\Bigl(\|u\|^2_{L^\infty}\,
\|\nabla\vrho\|^2_{L^\infty}
+\|\nabla\varrho\|_{L^\infty}^4+\|\nabla^2\varrho\|_{L^\infty}^2 \Bigr)
\,\|\vrho\|_{B^s_{p,r}}\,d\tau.
\end{align*}

By decomposing $\lambda$ into $\lambda(1)$ and $\lambda-\lambda(1)$,  one has
$$
\left\|\lambda(\rho)\,\nabla\pi\right\|_{\wtilde{L}^1_t(B^s_{p,r})}   \leq C\,\left(\left\|\nabla\pi\right\|_{\wtilde{L}^1_t(B^s_{p,r})}\,+\,\int^t_0\|\nabla\pi\|_{L^\infty}\,
\|\vrho\|_{B^s_{p,r}}\,d\tau\right)\,.
$$
One then uses \eqref{linest:pi,origin} to bound $\|\nabla\pi \|_{\wtilde{L}^1_t(B^s_{p,r})}$:   Lemma \ref{l:prod} help to bound the nonlinear term $\div((u+\nabla b)\cdot\nabla u)\equiv (\nabla u+\nabla^2 b):\nabla u$ by
\begin{eqnarray*}
  C\Bigl(\int^t_0\left(\|\nabla u\|_{L^\infty}+\left\|\nabla^2\vrho\right\|_{L^\infty}\right)\|u\|_{B^s_{p,r}}d\tau
  \,+ \,\int^t_0\|\nabla u\|^2_{L^\infty}\,\|\vrho\|_{B^s_{p,r}}\,d\tau
  \,+\,\eps \|\vrho\|_{\wtilde{L}^1_t(B^{s+2}_{p,r})}\Bigr)\,.
\end{eqnarray*}
Finally, Lemma \ref{l:tilde}   entails the control for the commutator term:
\begin{eqnarray*}
\left\|2^{j(s-1)}\int^t_0\|\div([\lambda,\Delta_j]\nabla\pi)\|_{L^p}d\tau\right\|_{\ell^r}
\lesssim\,\int^t_0
\left(\left\|\nabla\vrho\right\|_{L^\infty}\,\left\|\nabla\pi\right\|_{B^{s-1}_{p,r}}\,+\,
\left\|\nabla\vrho\right\|_{B^{s-1}_{p,r}}\,\left\|\nabla\pi\right\|_{L^\infty}\right)d\tau .
\end{eqnarray*}
Hence, interpolation inequality for $\|\nabla\pi\|_{B^{s-1}_{p,r}}$ between $B^{-\sigma}_{p,\infty}$ and $B^s_{p,r}$ helps us to   gather
\begin{eqnarray}
 \left\|\nabla\pi\right\|_{\wtilde{L}^1_t(B^s_{p,r}) }
 & \lesssim  & C(\|\nabla\varrho\|_{L^\infty_t(L^\infty)},s,\sigma)
     \|  \nabla\pi\|_{ {L}^1_t(B^{-\sigma}_{p,\infty})}   \label{est_cc:pi} \\
& & \;+\,\int^t_0\Bigl(\|\nabla u\|_{L^\infty}+\left\|\nabla^2\vrho\right\|_{L^\infty}+\|\nabla\vrho\|^2_{L^\infty}
    +\|\nabla\vrho\|^4_{L^\infty}\Bigr)
\|u\|_{B^s_{p,r}}d\tau\,\nonumber \\
& &  +\,\int^t_0\Bigl(\|\nabla u\|^2_{L^\infty}+\|u\|^2_{L^\infty}\|\nabla\vrho\|^2_{L^\infty}+\|\nabla^2\vrho\|^2_{L^\infty}+\|\nabla\pi \|_{L^\infty}\Bigr)
\|\vrho\|_{B^s_{p,r}}d\tau \nonumber \\
& & \qquad+\,\Bigl(1+\|\nabla\vrho\|_{L^\infty_t(L^\infty)}+\|u\|_{L^\infty_t(L^\infty)}\Bigr)
\|\vrho\|_{\wtilde{L}^1_t(B^{s+2}_{p,r})}\,. \nonumber
\end{eqnarray}

In the end, we discover from \eqref{est:til_u_rough} that $\|u\|_{\wtilde{L}^\infty_t(B^s_{p,r})}$ satisfies also
Inequality  \eqref{est_cc:pi}, just with an additional term $\|u_0\|_{B^s_{p,r}}$  on the right-hand side.
Recalling Estimate \eqref{est_cc:vrho} for the density,  we can replace $\|\vrho\|_{\wtilde{L}^1_t(B^{s+2}_{p,r})}$ in
Inequality \eqref{est_cc:pi} by the right-hand side of it.

Thus,  we can sum up \eqref{est_cc:vrho} and the (modified)
estimate \eqref{est_cc:pi} for the velocity $u$, yielding Lemma \ref{l:cc} by Gronwall's Lemma.


\subsection{Lower bounds for the lifespan of the solution} \label{s:life}

The aim of the present subsection is analyzing the lifespan of the solutions to system \eqref{system}.
We want to show, as carefully as possible, the dependence of the lifespan $T$ on the \textit{initial data}.

We will use freely
the following inequalities:
\begin{equation}\label{lifeprod}
\|ab\|_{B^{s-i}_{p,r}}\lesssim \|a\|_{B^{s-i}_{p,r}} \|b\|_{B^{s-i}_{p,r}},\quad
\|a^2\|_{B^s_{p,r}}\lesssim \|a\|_{L^\infty} \|a\|_{B^s_{p,r}},\quad
\|a\|_{L^\infty},\|\nabla a\|_{L^\infty}\lesssim \|a\|_{B^s_{p,r}},
\quad i=0,1,
\end{equation}
  and thanks to Conditions \eqref{index:s,r} and \eqref{index:p},
\begin{equation}\label{lifeineq:L2}
\|ab\|_{L^2}\leq \|a\|_{L^4} \|b\|_{L^4}\lesssim \|a\|_{B^{s-1}_{p,r}} \|b\|_{B^{s-1}_{p,r}}.
\end{equation}
By embedding results, without any loss of generality, throughout this subsection we will always assume $(s,p,r)=(1+d/4,4,1)$.
For notation convenience, we define $R_0\,:=\,\|\vrho_0\|_{B^{1+d/4}_{4,1}}$ and $U_0\,:=\,\|u_0\|_{B^{1+d/4}_{4,1}}$,
$$
R(t)\:=\,\left\|\vrho\right\|_{L^\infty_t(B^{1+d/4}_{4,1})}\,,\qquad S(t)\,:=\,\left\|\vrho\right\|_{L^1_t(B^{3+d/4}_{4,1})}
\qquad\mbox{ and }\qquad U(t)\,:=\,\left\|u\right\|_{L^\infty_t(B^{1+d/4}_{4,1})}\,.
$$
Then $S'(t)\equiv \|\varrho(t)\|_{B^{s+2}_{p,1}}$  controls high regularity of the density.

From \eqref{est_cc:vrho}, we infer that
$$
 R(t)\,+\,S(t) \,\leq\,C\Bigl(R_0\,+\,\int^t_0R  (U +1 )\,d\tau\,+\,
\int^t_0R^3 \,d\tau\Bigr).
$$
Now, if we define
\begin{equation} \label{def_life:T_R}
 T_R\;:=\;\sup\Bigl\{t\,>\,0\;\;\bigl|\;\;\int^t_0R^3(\tau)\,d\tau\;\leq\;2\,R_0\Bigr \}\,,
\end{equation}
by Gronwall's Lemma, for $t\in [0,T_R]$ and for some large enough $C$, we get
\begin{equation} \label{est_life:R-S}
R(t)\,+\,S(t)\,\leq\,C\,R_0\,\mc{E}(t),
\hbox{ with } \mc{E}(t):=\exp\left(C\,\int^t_0\bigl(1+U(\tau)\bigr)\,d\tau\right)\,.
\end{equation}


Next we aim to bound $U$, just as in the last subsection \ref{s:conti}.
Firstly, one applies the classical commutator estimates to the commutator $\cR_j=[u+\nabla b,\dj]\cdot\nabla u$ to arrive at
\begin{equation} \label{est_life:u_rough}
 U(t)\,\leq\,C\biggl(U_0
 \,+\,
\int^t_0 U^2  \,+\, U \,\|\varrho\|_{B^{s+1}_{p,1}} \,d\tau
 \,+\,\int^t_0\|h\|_{B^s_{p,1}}\,d\tau\,
 +\,\int^t_0\left\|\lambda\,\nabla\pi\right\|_{B^s_{p,1}}\,d\tau\biggr).
\end{equation}
Next, since the interpolation inequality for the  embeddings
$B^s_{p,1}\hra B^{s-1}_{p,1}\hra L^2$ holds, from estimate \eqref{linest:pi,origin} for $\pi$ one derives that, for some $\delta>1$,
$$
\|\nabla\pi\|_{B^s_{p,1}}\,\leq\,C\left((1+R^\delta)\left\|\nabla\pi\right\|_{L^2}
\,+\,(1+R)
\,\|(\nabla h, \nabla(u+\nabla b) : \nabla u)\|_{B^{s-1}_{p,1}}\right).
$$
Then, estimates \eqref{lifeprod} and \eqref{lifeineq:L2} help to    bound $h$, $\nabla(u+\nabla b):\nabla u$ and $\nabla\pi$ as follows:
$$\displaylines{
\|h\|_{B^s_{p,1}},\, \|\nabla(u+\nabla b):\nabla u\|_{B^{s-1}_{p,1}}\,\lesssim\,U\|\vrho\|_{B^{s+1}_{p,1}}+UR\|\vrho\|_{B^{s+1}_{p,1}}+
R^2S'+RS'+U^2,
\cr
 \|\nabla \pi\|_{L^2}
 \lesssim \|(h, \nabla(u+\nabla b):\nabla u)\|_{L^2}
 \lesssim UR^2+U\|\vrho\|_{B^{s+1}_{p,1}}+U^2+R^3+RS'+UR.
}$$
Now we use interpolation to write $\|\vrho\|_{B^{s+1}_{p,1}}\lesssim R^{1/2}(S')^{1/2}$, and Young inequality to separate
the term $S'$. Hence, from \eqref{est_life:u_rough} and previous inequalities one infers that
\begin{equation*}
U(t)\,\lesssim U_0+\int^t_0\left(1+R^{\delta+1}\right)\biggl(U(R^2+R)+U^2(1+R^3)+R^3+S'(1+R^2)\biggr)d\tau.
\end{equation*}

Let us restrict now to the interval $[0,T_R]$:   from \eqref{est_life:R-S} and the
previous estimate for $U$ we get (possibly taking a bigger $C$)
\begin{align*}
U(t)& \lesssim \mc{E}(t)\left(U_0+\left(1+R_0^{\delta+4}\right)\int^t_0\bigl(U+U^2+S'+R^3\bigr)d\tau\right) \\
& \lesssim \mc{E}(t)\left(1\,+\,U_0\,+\,R_0^{\delta+5}\,+\,\left(1+R_0^{\delta+4}\right)\int^t_0\bigl(U+U^2\bigr)d\tau\right)\,.
\end{align*}
Now we define
\begin{equation} \label{def_life:T_U}
T_U\;:=\;\sup \Bigl\{t>0\;\bigl|\;
\mc{E}(t)\leq 2,\,
 \left(1+R_0^{\delta+4}\right)\int^t_0\left(U+U^2\right)\,d\tau\,
 \leq\,2\left(1+U_0+R_0^{\delta+4}\right)\,
\; \Bigr\}.
\end{equation}
Then, for $0\leq t\leq \min\{T_R, T_U\}$, we infer the estimate
$$
U(t)\,\leq C\left(1\,+\,U_0\,+\,R_0^{\delta+5}\right)\,.
$$

Finally, by a bootstrap argument, it is easy to check that the time $T$, defined by \eqref{est:lifespan} in Theorem \ref{th:lifespan},  with sufficiently small $L$, is less than $T_R$ and $T_U$.
Hence Theorem \ref{th:lifespan} follows.

\appendix
\section{ Proof of Lemma \ref{l:tilde} }
\setcounter{equation}{0}
\numberwithin{equation}{section}
In the appendix we will prove Estimate \eqref{conest:comm-deri}.
The following classical properties will be used freely throughout this section:
\begin{itemize}
\item for any $u\in\cS',$ the equality $u=\sum_{j}\dj u$ holds true in $\cS'$;
\item for all $u$ and $v$ in $\cS',$
the sequence
$(S_{j-1}u\,\dj v)_{j\in\N}$ is spectrally
supported in dyadic annuli.
\end{itemize}

First of all, let us recall an easy version of Young inequality:
\begin{equation} \label{est:Young}
 a\,b\;\;\leq\;\;\theta\,\,\veps^{-(1-\theta)/\theta}\,\,a^{1/\theta}\,+\,(1-\theta)\,\,\veps\,\,b^{1/(1-\theta)},
 \quad
 \forall \theta\in\,]0,1[\,, \veps>0, \, a,b\in \R^+.
\end{equation}
 We decompose the commutator by use of Bony's paraproduct:
\begin{equation} \label{eq:comm_dec}
[\vphi,\Delta_j]\cdot\nabla\psi\,=\,R^1_j(\vphi,\psi)\,+\,R^2_j(\vphi,\psi)\,+\,R^3_j(\vphi,\psi)\,+\,R^4_j(\vphi,\psi)\,+\,
R^5_j(\vphi,\psi)\,,
\end{equation}
where, setting $\wtilde{\vphi}=(\Id-\Delta_{-1})\vphi$, we have defined
\begin{eqnarray*}
 R^1_j(\vphi,\psi) & := & \left[T_{\wtilde{\vphi}},\Delta_j\right]\cdot\nabla\psi \\
R^2_j(\vphi,\psi) & := & T'_{\Delta_j\nabla\psi}\wtilde{\vphi}\;=\;\sum_k S_{k+2}\dj\nabla\psi\,\cdot\Delta_k\wtilde{\vphi} \\
R^3_j(\vphi,\psi) & := & -\Delta_jT_{\nabla\psi}\wtilde{\vphi} \\
R^4_j(\vphi,\psi) & := & -\Delta_jR(\wtilde{\vphi},\nabla\psi) \\
R^5_j(\vphi,\psi) & := & \left[\Delta_{-1}\vphi,\Delta_j\right]\cdot\nabla\psi\,.
\end{eqnarray*}

One finds easily
$$
R^1_j(\vphi,\psi)\,=\,\sum_{|\nu-j|\leq 1}\int_{\R^d_z}2^{-j}\left(\int^1_0h(z)\,z\,
\cdot\nabla S_{\nu-1}\wtilde{\vphi}(x-2^{-j}\lambda z)\,\cdot\Delta_\nu\nabla\psi(x-2^{-j}z)\,d\lambda\right)dz\,.
$$
This ensures that
$$
\mc R^1:=\Bigl\| \Bigl( 2^{js}\int^t_0 \|\nabla R^1_j\|_{L^p}\Bigr)_j\Bigr\|_{\ell^r}
\,\lesssim\,\Bigl\| \Bigl(  2^{js}\int^t_0 \|\nabla S_{j-1}\wtilde{\vphi} \|_{L^\infty}\,
 \|\Delta_j\nabla\psi \|_{L^p}\,d\tau\Bigr)_j\Bigr\|_{\ell^r}.
$$
We apply Young inequality \eqref{est:Young} to the integrand on the right hand side to get, for some constant $C$:
\begin{equation} \label{est:tilde_1}
 \mc R^1\,\leq\,\frac{C\,\theta}{\veps^{(1-\theta)/\theta}}\,
\int^t_0\|\nabla\vphi\|^{1/\theta}_{L^\infty}\,\|\psi\|_{B^{s_1}_{p,r}}\,d\tau\,+\,\,(1-\theta)\,\veps\,
\|\psi\|_{\wtilde{L}^1_t(B^{s_2}_{p,r})}\,.
\end{equation}

Let us now handle
\begin{align*}
\mc R^2
&:= \Bigl\|2^{j\,s}\int^t_0 \|\nabla R^2_j \|_{L^p}\,d\tau\Bigr\|_{\ell^r} \\
& \lesssim
\Bigl\|2^{j\,s}\int^t_0\sum_{\mu\geq j-2}
\Bigl(\|\nabla^2 S_{\mu+2} \Delta_j\psi\|_{L^\infty}\,
\|\Delta_\mu\wtilde\vphi\|_{L^p}
+ \|S_{\mu+2}\nabla\Delta_j\psi\|_{L^\infty}\,
\|\nabla\Delta_\mu\wtilde\vphi\|_{L^p} \Bigr) \,d\tau\Bigr\|_{\ell^r} \\
&\lesssim
\Bigl\|2^{j\,s}\int^t_0\sum_{\mu\geq j-2}\|\nabla\Delta_j\psi\|_{L^\infty}\,
\|\nabla\Delta_\mu\wtilde\vphi\|_{L^p}\,d\tau\Bigr\|_{\ell^r} \\
&\lesssim
\Bigl\|\int^t_0\|\nabla\psi\|_{L^\infty}\,
\sum_{\mu\geq j-2}2^{(j-\mu)s}2^{\mu s}
\|\nabla\Delta_\mu\wtilde\vphi\|_{L^p}\,d\tau\Bigr\|_{\ell^r}\,.
\end{align*}
We just do exactly as above (the way to obtain \eqref{est:tilde_1} ): if $s>0$, then we have
\begin{equation} \label{est:tilde_2}
 \mc R^2\,\leq\,\frac{C\,\eta}{\veps^{(1-\eta)/\eta}}\,
\int^t_0\|\nabla\psi\|^{1/\eta}_{L^\infty}\,\|\nabla\tilde\varphi\|_{B^{\sigma_1-1}_{p,r}}\,d\tau
\,+\,\,(1-\eta)\,\veps\,\|\nabla\tilde\varphi\|_{\wtilde{L}^1_t(B^{\sigma_2-1}_{p,r})}\,.
\end{equation}

Moreover, since
$$
\mc R^3\,:=\,\Bigl\|2^{j\,s}\int^t_0\left\|\nabla R^3_j\right\|_{L^p}\,d\tau\Bigr\|_{\ell^r}\,\leq\,
\Bigl\|2^{j(s+1)}\int^t_0\sum_{\mu\sim j}\|S_{\mu-1}\nabla\psi\|_{L^\infty}\,
\|\Delta_\mu\wtilde\vphi\|_{L^p}\Bigr\|_{\ell^r}
$$
we can immediately see that \eqref{est:tilde_2} holds also for $\cR^3$.
Similarly, \eqref{est:tilde_2}   follows immediately for
$$
\mc R^4
\,:=\,\Bigl\|2^{j\,s}\int^t_0\left\|\nabla R^4_j\right\|_{L^p}\,d\tau\Bigr\|_{\ell^r}
\lesssim
\,\Bigl\|\int^t_0 \|\nabla\psi\|_{L^\infty}\,\sum_{\mu\geq j-2}2^{(j-\mu)(s+1)}\left(2^{\mu(s+1)}\|\Delta_\mu\wtilde\vphi\|_{L^p}\right)\,
\,d\tau\Bigr\|_{\ell^r}.
$$
Finally, the last term
$$
\mc R^5\,:=\,\Bigl\|2^{j\,s}\int^t_0\left\|\nabla R^5_j\right\|_{L^p}\,d\tau\Bigr\|_{l^r}
$$
can be handled as $\mc R^1$, leading us to the same estimate as \eqref{est:tilde_1} and so to the end of the proof.

\subsubsection*{Acknowledgements}
The main part of the work was prepared when the first author was a post-doc at BCAM - Basque Center for Applied Mathematics, and
the second author was a Ph.D. student at LAMA - Laboratoire d'Analyse et de Math\'ematiques Appliqu\'ees, UMR 8050,
Universit\'e Paris-Est. They want to aknowledge both these institutions.

The first author was partially supported by Grant MTM2011-29306-C02-00, MICINN, Spain,
ERC Advanced Grant FP7-246775 NUMERIWAVES, ESF Research Networking Programme OPTPDE and Grant PI2010-04 of the Basque Government.
During the last part of  the work, he was also supported by the project ``Instabilities in Hydrodynamics'',
funded by the Paris city hall (program ``\'Emergences'') and the
Fondation Sciences Math\'ematiques de Paris.

The second author  was partially supported by the  project ERC-CZ
LL1202, funded by the Ministry of Education, Youth and Sports of the Czech Republic.



The first author is member of the Gruppo Nazionale per l'Analisi Matematica, la Probabilit\`a e le loro Applicazioni (GNAMPA) of the
Istituto Nazionale di Alta Matematica (INdAM).

\bigbreak\bigbreak

\vspace{.5cm}

\begin{flushleft}
\textsl{Francesco Fanelli} \\ \vspace{.1cm} 
{\small \textit{Institut de Math\'ematiques de Jussieu-Paris Rive Gauche -- UMR 7586} \\
 \textsc{Universit\'e Paris-Diderot -- Paris 7} \\
  B\^atiment Sophie-Germain, case 7012 \\
   56-58, Avenue de France \\
   75205 Paris Cedex 13 -- FRANCE \\ \vspace{.1cm}
   E-mail: \texttt{fanelli@math.jussieu.fr} }

 \bigbreak

\textsl{Xian Liao} \\ \vspace{.1cm} 
{\small \textit{Academy of Mathematics \& Systems Science} \\
 \textsc{Chinese Academy of Sciences}\\
  55 Zhongguancun East Road  \\
  100190 Beijing -- P.R. CHINA \\ \vspace{.1cm}
E-mail: \texttt{xian.liao@amss.ac.cn} }
\end{flushleft}

\end{document}